\documentclass[10pt,a4paper]{amsart}
\usepackage{latexsym,amssymb,amsmath,amsthm,amsfonts,enumerate,verbatim,xspace,exscale}
\usepackage[mathscr]{eucal}

\input xy
\xyoption{all}
\CompileMatrices
\UseComputerModernTips

\theoremstyle{plain}
\newtheorem{theorem}{Theorem}[section]
\newtheorem{lemma}[theorem]{Lemma}
\newtheorem{proposition}[theorem]{Proposition}
\newtheorem{corollary}[theorem]{Corollary}

\theoremstyle{definition}

\newtheorem{definition}[theorem]{Definition}

\newtheorem{rem}[theorem]{Remark}

\newtheorem*{name}{\theoremname}
\newcommand{\theoremname}{testing}

\newtheorem{stepp}{Remark}[section]

\newtheorem{exx}[theorem]{Example}

\allowdisplaybreaks

\newcommand{\caps}[1]{\textup{\textsc{#1}}}

\providecommand{\bysame}{\makebox[3em]{\hrulefill}\thinspace}

\newcommand{\iso}{isomorphism\xspace}

\newcommand{\diffeo}{diffeomorphism\xspace}

\newcommand{\bb}{\mathbb}
\newcommand{\ov}[1]{\mbox{$\overline{#1}$}}
\newcommand{\up}{\upshape}

\newcommand{\x}{$\hfill\Box$}

\providecommand{\xrighto}[1]{\mbox{$\;\xrightarrow{#1}\;$}}

\newcommand{\longto}{\longrightarrow}
\newcommand{\hookto}{\hookrightarrow}
\newcommand{\toto}{\twoheadrightarrow}

\def\vv<#1>{\langle#1\rangle}

\newcommand{\tr}{\mbox{$\textup{Tr}$}}
\newcommand{\diag}[1]{\mbox{$\textup{diag}(#1)$}}

\newcommand{\im}{\mbox{$\text{\up{im}}\,$}}

\newcommand{\ddim}{\mbox{$\text{ddim}\,$}}

\providecommand{\sign}{\mbox{$\text{\up{sign}}\,$}}

\newcommand{\dd}[2]{\mbox{$\frac{\partial #2}{\partial #1}$}}

\providecommand{\del}{\partial}

\newcommand{\om}{\omega}
\newcommand{\Om}{\Omega}
\newcommand{\var}{\varphi}

\newcommand{\lam}{\lambda}
\newcommand{\Lam}{\Lambda}

\newcommand{\wt}[1]{\mbox{$\widetilde{#1}$}}

\newcommand{\bsc}{\mbox{$\bigsqcup$}}

\newcommand{\R}{\mbox{$\bb{R}$}}
\newcommand{\C}{\mbox{$\bb{C}$}}

\newcommand{\by}[2]{\mbox{$\frac{#1}{#2}$}}
\newcommand{\cinf}{\mbox{$C^{\infty}$}}
\providecommand{\set}[1]{\mbox{$\{#1\}$}}
\newcommand{\subeq}{\subseteq}

\newcommand{\X}{\mathfrak{X}}

\newcommand{\curv}{\mbox{$\textup{Curv}$}}

\newcommand{\mf}{manifold\xspace}

\newcommand{\ver}{\mbox{$\textup{Ver}$}}
\newcommand{\hor}{\mbox{$\textup{Hor}$}}

\newcommand{\ann}{\mbox{$\textup{Ann}\,$}}

\newcommand{\hcm}{\mbox{$H_{\textup{CM}}$}} 
\newcommand{\hcml}{\mbox{$H_{\textup{CM}}^{(L)}$}} 
\newcommand{\hfree}{\mbox{$H_{\textup{free}}$}}
\newcommand{\bra}[1]{\mbox{$\{#1\}$}}
\newcommand{\ham}[1]{\mbox{$\smash{\nabla^{\om}_{#1}}$}}

\newcommand{\spr}[1]{\mbox{$/\negmedspace/_{#1}$}}
\newcommand{\sporb}{\mbox{$/\negmedspace/_{\mathcal{O}}$}}
\newcommand{\momap}{momentum map\xspace}
\newcommand{\symplecto}{symplectomorphism\xspace}
\newcommand{\omkks}{\mbox{$\Om^{\mathcal{O}}$}}

\newcommand{\csa}{Cartan subalgebra\xspace}

\newcommand{\lie}[1]{\mbox{$\text{\up{Lie}}(#1)$}}

\newcommand{\gu}{\mathfrak{g}}
\newcommand{\ho}{\mathfrak{h}}

\newcommand{\ko}{\mathfrak{k}}
\newcommand{\lo}{\mathfrak{l}}
\newcommand{\mo}{\mathfrak{m}}

\newcommand{\po}{\mathfrak{p}}

\newcommand{\too}{\mathfrak{t}}

\newcommand{\conj}{\mbox{$\text{\up{conj}}$}}
\newcommand{\Ad}{\mbox{$\text{\upshape{Ad}}$}}
\newcommand{\ad}{\mbox{$\text{\upshape{ad}}$}}

\newcommand{\orb}{\mbox{$\mathcal{O}$}}

\newcommand{\SL}{\mbox{$\textup{SL}$}}

\newcommand{\SU}{\mbox{$\textup{SU}$}}
\newcommand{\su}{\mbox{$\mathfrak{su}$}}
\newcommand{\SO}{\mbox{$\textup{SO}$}}

\newcommand{\FO}{\mbox{$\mathcal{F}$}}

\newcommand{\ilat}{\mbox{$\mathcal{IL}$}}

\newcommand{\WW}{\mbox{$\mathcal{W}$}}

\newcommand{\lorb}[1]{\mbox{${#1}_{(L)}^{\mathcal{O}}$}}

\newcommand{\bscorb}{\mbox{$\bsc_{q\in Q}\mathcal{O}\cap\ann\gu_{q}$}}
\newcommand{\ine}{\mbox{$\mathbb{I}$}}

\newcommand{\kks}{Kazdhan, Kostant, Sternberg\xspace}
\newcommand{\aklm}{Alekseevsky, Kriegl, Losik, Michor\xspace}

\title[Cotangent Bundle Reduction \& Calogero-Moser Systems]{Singular Cotangent Bundle  Reduction\\
       \&\\
       Spin Calogero-Moser Systems}
\author{Simon Hochgerner}
\address{Fakult\"at f\"{u}r Mathematik\\
  Universit\"{a}t Wien\\
  Nordbergstrasse~15\\
  A-1090 Vienna, Austria}
\email{simon.hochgerner@univie.ac.at}  
\urladdr{http://www.mat.univie.ac.at/\~{}simon}
\thanks{This work is supported by Fonds zur F\"{o}rderung der
       wissenschaftlichen Forschung, Projekt P~14195-MAT}
\keywords{Cotangent bundle reduction, singular reduction, Hamiltonian
  systems, Calogero-Moser systems}
\subjclass[2000]{53D20}

\begin{document}

\begin{abstract}
We develop a bundle picture
for singular symplectic quotients of cotangent bundles acted upon by
cotangent lifted actions 
for the case that the
configuration manifold is of single orbit type.
Furthermore,
we give a
formula for the reduced symplectic form in this setting. 
As an application of this bundle picture we
consider Calogero-Moser systems with spin associated to polar representations
of compact Lie groups.
\end{abstract}
\maketitle

\section{Introduction}

Let $G$ be a Lie group acting properly on a configuration manifold
$Q$. 
Consider the cotangent lifted $G$-action on $T^*Q$.
This action is
Hamiltonian 
with respect to the 
standard exact symplectic form on $T^*Q$ and  
with equivariant \momap 
denoted by $\mu: T^*Q\to\gu^*$.
Assume $\orb$ is a
coadjoint orbit contained in the image of $\mu$. 
The first part of this paper 
is concerned with the study of the singular symplectic
quotient
\[
 \mu^{-1}(\orb)/G =: T^{*}Q\sporb G.
\]
Indeed, this quotient cannot be a smooth manifold, in
general, since we do not assume the $G$-action on $Q$ to be free. 
However, we can apply the theory of 
singular symplectic
reduction as developed by Sjamaar and Lerman \cite{SL91}, Bates and
Lerman \cite{BL97}, and Ortega and Ratiu \cite{OR04} (see also Theorem
\ref{thm:sing_spr}), and this exhibits $T^*Q\sporb G$ to be a Whitney
stratified space with strata of the form 
\[
 (\mu^{-1}(\orb)\cap(T^{*}Q)_{(L)})/G
 =:
 (T^{*}Q\sporb G)_{(L)}
\]
where $(L)$ is an element of the isotropy lattice of the $G$-action on
$T^{*}Q$.

One of the aims of this paper is to develop a bundle picture for the
reduced phase space $T^*Q\sporb G$, i.e., to obtain a fiber bundle
$T^*Q\sporb G\to T^*(Q/G)$ in a suitable (singular) sense. One can
hope to construct such a fiber bundle in the presence of an
additionally chosen generalized connection form (see
Section~\ref{sec:gauged_red}) on $Q\toto Q/G$. 
However, for reasons explained
in Remark~\ref{rem:stratQ} there cannot exist a surjection 
$T^*Q\sporb G\to T^*(Q/G)$ with locally constant fiber type
for general proper $G$-actions.
In Remark~\ref{rem:stratQ}(\ref{equ:F}) we state a weak substitute of
a bundle picture for $T^*Q\sporb G$.

To obtain a bundle picture and a useful description of the reduced
phase space $T^*Q\sporb G$
we have to assume that the base manifold $Q$  is of
single orbit type, that is, $Q=Q_{(H)}$ for a subgroup $H$ of
$G$. Assuming this we get a first result that says that, locally, 
\[
\xymatrix{
 {\orb\spr{0}H}
 \ar @{^{(}->}[r]&
 {T^{*}Q\sporb G}
 \ar[r]& 
 {T^{*}(Q/G)}
}
\]
is a symplectic fiber bundle (Theorem~\ref{thm:bun_pic}). 
This result is obtained by applying the Palais
Slice Theorem to the $G$-action on the base space $Q$, and then using the
Singular Commuting Reduction Theorem of Section~\ref{s:sg_com_red}. 
This is an inroad that was also taken by
Schmah~\cite{Sch04} to get a local description of $T^{*}Q\sporb G$. 

However, one can also give a global symplectic description of the
reduced phase space and this is done in Section \ref{sec:gauged_red}. 
This follows an approach that is generally called Weinstein
construction (\cite{Wei78}). 
(We do not consider the closely related construction of Sternberg~\cite{Ste77}.)
In the case that the $G$-action on
the configuration space $Q$ is free this global description was first given
by Marsden and Perlmutter~\cite{MP00}. Their result says that 
the choice of a principal bundle connection on $Q\toto Q/G$ yields a 
realization of
the symplectic
quotient $T^{*}Q\sporb G$ as a fibered product
\[
 T^{*}(Q/G)\times_{Q/G}(Q\times_{G}\orb),
\]
and they compute the reduced symplectic structure in terms of data
intrinsic to this realization -- \cite[Theorem 4.3]{MP00}. 

In the presence of a single non-trivial isotropy $(H)$ 
on the configuration space one
obtains a non-trivial isotropy lattice on $T^{*}Q$ 
whence the symplectic reduction of $T^*Q$ is to be carried out in the
singular context of \cite{SL91,BL97,OR04}.
The result is the following:
The choice of a generalized connection form 
(see Section~\ref{sec:gauged_red}) 
on $Q\toto Q/G$ yields a realization of 
each symplectic stratum $(T^{*}Q\sporb G)_{(L)}$ of the reduced space
as a fibered product
\[ 
  (\WW\sporb G)_{(L)}
  =
  T^{*}(Q/G)\times_{Q/G}(\bscorb)_{(L)}/G
\] 
where 
\[
 \WW :=
 (Q\times_{Q/G}T^{*}(Q/G))\times_{Q}\bsc_{q\in Q}\ann{\gu_{q}}
 \cong
 T^{*}Q
\]
as symplectic manifolds with a Hamiltonian $G$-action.
Moreover, we compute the reduced symplectic structure in terms
intrinsic to this realization. This is the content of 
Theorem~\ref{thm:WWred}. 
Even though we have to restrict to single orbit type
manifolds this result is quite general 
in the sense that it is valid for arbitrary
coadjoint orbits $\orb$. 

The first to have studied symplectic reduction of cotangent bundles
for non-free actions seems to have been Montgomery~\cite{Mon83}. 
Using the point reduction approach this paper gives conditions under
which the reduced phase space carries a smooth manifold structure. 
The
first to study this subject 
in the context of singular symplectic reduction as developed by \cite{SL91,BL97,OR04}
is Schmah~\cite{Sch04} 
who proves a cotangent bundle specific slice theorem at points whose
momentum values are fully isotropic.
The other
important paper on singular cotangent bundle reduction is by
Perlmutter, Rodriguez-Olmos and Sousa-Diaz~\cite{PRS03}. 
By restricting to do
reduction at fully isotropic values of the \momap $\mu:
T^{*}Q\to\gu^{*}$ they are able to drop all assumptions on the
isotropy lattice of the $G$-action on $Q$, and give a very complete
description of the reduced symplectic space. 

As an application of the bundle picture found in
Theorem~\ref{thm:WWred} we
consider Calogero-Moser systems with spin in Section \ref{sec:cms}.
In fact, it was an idea of \aklm \cite{AKLM03} to consider polar
representations of compact Lie groups $G$ on a Euclidean vector 
space $V$ to obtain new versions of
Calogero-Moser models. We make these ideas precise by using the
singular cotangent bundle reduction machinery. Thus let $\Sigma$ be a
section for the $G$-action in $V$, let $C$ be a Weyl chamber in
this section, and put $M:=Z_{G}(\Sigma)$. Under a strong but not
impossible condition on a chosen coadjoint orbit in $\gu^{*}$ we
get 
\[
 T^{*}V\sporb G
 =
 T^{*}C_{r}\times\orb\spr{0}M
\]
from the general theory (Theorem~\ref{thm:WWred}), 
where $C_{r}$ denotes the sub-manifold of
regular elements in $C$.  
This is the effective phase space of the spin
Calogero-Moser system. The corresponding Calogero-Moser function is obtained as a
reduced Hamiltonian from the free Hamiltonian on $T^{*}V$. The
resulting formula is  
\[
 \hcm(q,p,[Z])
 =
 \by{1}{2}\sum_{i=1}^{l}p_{i}^{2}
  + 
 \by{1}{2}\sum_{\lam\in R}
  \by{\sum_{i=1}^{k_{\lam}}z_{\lam}^{i}z_{\lam}^{i}}{\lam(q)^{2}}.
\]
This is made precise with the necessary notation in Section
\ref{sec:cms}.

It was first observed by Kazhdan, Kostant and Sternberg
\cite{KKS78} that one can obtain
Calogero-Moser models via Hamiltonian reduction of $T^*\gu$ or
$T^*G$. 
(See also Subsections~\ref{sub:cms-constr} and \ref{sub:kks}.) 
In fact, \cite{KKS78}
chose their example such that the reduction procedure yields a smooth
symplectic manifold, and moreover there appear no spin variables, i.e.,
$\orb\spr{0}H=\set{\textup{point}}$. 
This reduction approach was further pursued in \cite{AKLM03}
as alluded to above.
However, these spin models are Hamiltonian systems on non-smooth
spaces, and it was not clear in which sense one should  regard
the reduced systems as being Hamiltonian systems. 
Moreover, the precise form of the
singularities was not clear. These questions are answered by
Theorem~\ref{thm:WWred}: The reduced system is a stratified
Hamiltonian system, the strata of the reduced phase space are
described, and when the system is restricted to a stratum it is a
Hamiltonian system in the usual sense on this stratum. Moreover, the
stratification is Whitney whence by \cite{SL91}
the singularities which appear are
of conic form.

As an interesting side product (Remark~\ref{rem:r+A}) 
our (singular) cotangent bundle
reduction approach yields a connection to the $r$-matrix theoretic
construction of Calogero-Moser models of Li and
Xu~\cite{LX00,LX02}. This connection is new and
interesting since it explains in geometric terms 
why solutions of the classical dynamical
Yang-Baxter equation lead to Calogero-Moser systems.

Finally, we use a result on non-commutative integrability from Zung
\cite[Theorem 2.3]{Zung02} to show that these Calogero-Moser systems are
integrable in the non-commutative sense.

\textbf{Thanks.}
This paper is part of my PhD thesis written under
the supervision of Peter Michor. I am grateful to him for introducing me
to symplectic geometry and proposing the subject of Calogero-Moser
systems associated to polar representations of compact Lie
groups. Further, I wish to thank the Centre Bernoulli for their hospitality
in September 2004. This stay contributed a lot toward the
finishing of this paper. I am also thankful to Stefan Haller and Armin
Rainer for helpful remarks and comments. 
Finally, I want to thank the referees for their detailed report.

\section{Preliminaries and notation}\label{sec:prel}

All manifolds to be considered are Hausdorff, para-compact, 
finite dimensional, and smooth in the $\cinf$-sense. 
We do not assume that manifolds are connected but allow for finitely
many connected components of varying (finite) dimension. Thus the
dimension of a manifold is only locally constant. 
A proper
$G$-space $M$ is a manifold $M$ acted upon properly by a Lie group
$G$. 
For proper $G$-spaces the 
Slice Theorem~(\cite{Pal61,PT88,DK99}) holds, 
and we shall make frequent use of this theorem. 

Let $M$ be a proper $G$-space. An \caps{isotropy class} $(H)$ of the
$G$-action on $M$
is a conjugacy class of an isotropy subgroup $G_x$ of a point 
$x\in M$, 
that is, $(H) = \set{gG_xg^{-1}: g\in G}$. 
If we want to explicit that $(H)$ is the conjugacy class of $H$ with
respect to $G$ we shall write $(H)^G$.
The 
\caps{isotropy lattice} $\ilat(M)$ is defined to be the lattice consisting of all
isotropy classes $(H)$ of the $G$-space $M$. 
For $(H)\in\ilat(M)$ the \caps{orbit type} sub-manifold is
\begin{align}
 M_{(H)} := 
 M_{(H)^G} :=
 \set{x\in M: G_x \textup{ is conjugate to } H \textup{ within } G},
\end{align}
the \caps{symmetry type} sub-manifold is 
$M_H := \set{x\in M: G_x = H}$, and the \caps{fixed point} sub-manifold is 
$M^H := \set{x\in M: H\subset G_x}$. 
More generally,
let $N$ be a closed topological 
subspace of $M$ and let $K$ be a compact subgroup of $G$ which acts
(continuously)
on $N$. 
We may view $M$ also as a proper $K$-space, and thus obtain an orbit type stratification of 
$M$
with respect to the $K$-action. 
Since $N$ is $K$-invariant this induces a decomposition
of $N$ according to orbit types. 
Let $N_{(L_0)} = N_{(L_0)^K} = N\cap M_{(L_0)^K}$ 
be such an orbit type stratum, that is $L_0 = K_x$ for some $x\in N$.
With regard to $L_0\subset K$
we will be concerned with the \caps{generalized isotropy class}
\begin{align}
 (L_0)_K^G := \set{L\subset K: \textup{ there is } g\in G \textup{ s.t.\ } gL_0g^{-1} = L} 
\end{align}
and the corresponding generalized orbit type space 
\begin{align}
 N_{(L_0)_K^G} := \set{z\in N: K_z\in(L_0)_K^G}\label{def:weirdo-type}.
\end{align}
Clearly, $N_{(L_0)_K^G}$ itself 
decomposes into orbit type strata with respect to the induced $K$-action.

Using the Slice Theorem one can show (see
\cite{DK99}) that the stratification of $M$ into orbit type
sub-manifolds forms a Whitney stratification. Likewise the
stratification of the topological space $M/G$ into strata of the type
$M_{(H)}/G$ where $(H)\in\ilat(M)$
forms a Whitney stratification. Thus $M/G$ becomes a stratified
space, i.e., a stratified space such that the stratification satisfies
the Whitney conditions. See \cite{Mat70,GorMac88,Pfl01,Pfl01a} for
more on stratified spaces. If $X$ and $Y$ are stratified spaces a
\caps{stratified map} $\phi: X\to Y$ is a continuous map which respects the
stratifications, that is, the pre-image under $\phi$ of every stratum of $Y$  
decomposes into a union of strata of $X$. 
For example, the orbit projection map $\pi: M\toto M/G$ is a stratified
map. 

If the action is written as $l: G\times M\to M$, $(k,x)\mapsto
l(k,x)=l_{k}(x)=l^{x}(k)=k.x$ 
we can tangent bundle lift it via
$k.(x,v):=(k.x,k.v):=Tl_{k}.(x,v)=(l_{k}(x),T_{x}l_{k}.v)$ 
for $(x,v)\in TM$ to an action on $TM$. 
As the
action consists of transformations by diffeomorphisms 
it may also be lifted to
the cotangent bundle. This is the cotangent lifted action which is
defined by
$k.(x,p):=(k.x,k.p):=T^{*}l_{k}.(x,p)=
(k.x,T_{k.x}^{*}l_{k^{-1}}.p)$ where $(x,p)\in T^{*}M$. 
Our notation for the fundamental vector field is 
$
 \zeta_{X}(x) 
 := \zeta(x)(X)
 := \dd{t}{}|_{0}l(\exp(+tX),x)
  = T_{e}l^{x}(X)
$
where $X\in\gu$.

\section{Singular symplectic reduction}\label{s:sg_com_red}
Let $(M,\om)$ be a connected symplectic \mf, and $G$ a Lie group that
acts on $(M,\om)$ in a proper and Hamiltonian fashion such that there is an
equivariant \momap $J: M\to\gu^{*}$.

The very strong machinery of singular symplectic reduction is 
(for the case of compact $G$) due to
Sjamaar and Lerman~\cite{SL91} who prove that the singular symplectic
quotient is a Whitney stratified space that has symplectic manifolds
as its strata. This result which is the Singular Reduction Theorem was
then generalized to the case of proper actions by Bates and
Lerman~\cite{BL97}, 
Ortega and Ratiu~\cite{OR04}, and others. 

\begin{theorem}[Singular symplectic reduction]\label{thm:sing_spr}
Let $(H)$ be in the isotropy lattice of the $G$-action on $M$, and
suppose that $J^{-1}(\orb)\cap M_{(H)}\neq\emptyset$ for a coadjoint
orbit $\orb\subeq\gu^{*}$. 
Then the following are true.
\begin{itemize}
\item
The subset $J^{-1}(\orb)\cap M_{(H)}$ is an initial sub-manifold of
$M$.
\item 
The topological quotient $(J^{-1}(\orb)\cap M_{(H)})/G$ has a unique
smooth structure such that the projection map
\[
\xymatrix{
 {J^{-1}(\orb)\cap M_{(H)}}
 \ar @{->>}[r]^-{{\pi}}&
 {(J^{-1}(\orb)\cap M_{(H)})/G}
}
\]
is a smooth surjective submersion. 
\item
Let $\iota: J^{-1}(\orb)\cap M_{(H)}\hookto M$ denote the inclusion
mapping. Then $(J^{-1}(\orb)\cap M_{(H)})/G$ carries a symplectic
structure $\om_{0}$ which is uniquely characterized by the formula
\[ 
 \pi^{*}\om_{0}  
 = 
  \iota^{*}\om
  - (J|(J^{-1}(\orb)\cap M_{(H)}))^{*}\Om^{\mathcal{O}}
\]
where $\Om^{\mathcal{O}}$ is the canonical (positive Kirillov-Kostant-Souriau)
symplectic form on $\orb$.
\item 
Consider a $G$-invariant function $H\in\cinf(M)^{G}$. Then the flow of
the Hamiltonian vector field $\ham{H}$ leaves the connected components
of $J^{-1}(\orb)\cap M_{(H)}$ invariant. Moreover, $H$ factors to a
smooth function $h$ on the quotient 
$(J^{-1}(\orb)\cap M_{(H)})/G$. 
Finally, $\ham{H}$ and the Hamiltonian vector field to $h$
are related via the canonical projection $\pi$, whence the flow of the former
projects to the flow of the latter.
\item 
The collection of all strata of the form $(J^{-1}(\orb)\cap M_{(H)})/G$
constitutes a Whitney stratification of the topological space
$J^{-1}(\orb)/G$.
\end{itemize}
\end{theorem}
 
\begin{proof}
This theorem is contained in \cite[Section 8]{OR04}. See also
\cite[Corollary 14]{BL97} and \cite{SL91}.
\end{proof}

In fact, \cite{OR04} state the above theorem only for connected
components of strata. This is so because they do not require
the \momap $J$ to be equivariant with respect to the co-coadjoint action
on $\gu^*$. 

As a matter of convention we write shorthand $M\sporb G :=
J^{-1}(\orb)/G$ for the reduced space of $M$ with respect to the
Hamiltonian action by $G$. If $\orb$ is the coadjoint orbit passing
through $\alpha$ then we shall also abbreviate
$J^{-1}(\alpha)/G_{\alpha} = M\spr{\alpha}G = M\sporb G$.

\begin{theorem}[Singular commuting reduction]\label{thm:sg_comm_red}
Let $G$ and $H$ be Lie groups that act properly and by symplectomorphisms
on $(M,\om)$ with momentum maps $J_{G}$ and $J_{H}$
respectively. Assume that the actions commute, that $J_{G}$ is
$H$-invariant, and that $J_{H}$ is $G$-invariant. Let $\alpha\in\gu^{*}$ be in
the image of $J_{G}$ and $\beta\in\ho^{*}$ in the image of $J_{H}$.

Then the $G$ action drops to a Poisson action on $M\spr{\beta}H$ and
$J_{G}$ factors to a \momap $j_{G}$ for the induced action. Likewise, 
the $H$ action drops to a Poisson action on $M\spr{\alpha}G$ and
$J_{H}$ factors to a \momap $j_{H}$ for the induced
action. Furthermore, we have
\[
 (M\spr{\alpha}G)\spr{\beta}H 
 \cong
 M\spr{(\alpha,\beta)}(G\times H)
 \cong
 (M\spr{\beta}H)\spr{\alpha}G
\]
as symplectic stratified spaces.
\end{theorem}

\begin{proof}
An outline of a proof of this result is given in \cite[Section~4]{SL91} for the case that 
$G$ and $H$ are compact. Using the machinery of singular symplectic reduction 
for proper Hamiltonian actions as described in \cite{OR04} the proof 
of \cite{SL91} extends to the more general setting. 
\end{proof}

\section{The bundle picture}\label{sec:bun-pic}

From now on let $G$ be a Lie group acting properly from the left on a
manifold $Q$. 
The $G$ action then induces a Hamiltonian action
on the cotangent bundle $T^{*}Q$ by cotangent lifts. This means that
the lifted action respects the canonical symplectic form
$\Om=-d\theta$ on $T^{*}Q$ where $\theta$ is the Liouville form on
$T^{*}Q$, and, 
moreover, there is an
equivariant \momap 
$\mu:
T^{*}Q\to\gu^{*}$ given by 
$\vv<\mu(q,p),X>
 = \theta(\zeta^{T^{*}Q}_{X})(q,p)
 = \vv<p,\zeta_{X}(q)>$ 
where
$(q,p)\in T^{*}Q$, $X\in\gu$, $\zeta_{X}$ is the fundamental
vector field associated to the $G$-action on $Q$, 
and $\zeta^{T^{*}Q}_{X}\in\X(T^{*}Q)$ is the
fundamental vector field associated to the cotangent lifted action. 

In this section we want to apply the Slice Theorem (\cite{PT88,DK99,OR04})
 to the action of $G$ on $Q$ to get a local model of
the singular symplectic reduced space $T^{*}Q\sporb G =
\mu^{-1}(\orb)/G$ where $\orb$ is a coadjoint orbit in the image of
$\mu$.  

Thus we consider a tube $U$ in $Q$ around an orbit $G.q$ with
$G_{q}=H$. And we denote the slice at $q$ by $S$ such that 
\[
 U\cong G\times_{H}S
\]
as $G$-spaces. Here the action on $U$ is given by the restriction of
the $G$-action on $Q$ to the invariant neighborhood $U$ of $G.q$.
On the other hand, the
action on $G\times_H S$ is given by $g.[(k,s)]_H = [(gk,s)]_H$ where
$g\in G$ and $[(k,s)]_H$ denotes the class of $(k,s)\in G\times S$ in
$G\times_H S$. Moreover, note that the $H$-action on $G\times S$ is
given by $h.(k,s) = (kh^{-1},h.s)$ where $h\in H$. 
(See \cite{PT88,DK99,OR04}.)
In particular, it follows that $U/G\cong S/H$ as stratified spaces
with smooth structure. (See \cite{Pfl01,Pfl01a}.)

Assume for a moment that the action by $G$ on $U$ is free, whence
$U\cong G\times S$. Let $\mu: T^{*}U\to\gu^{*}$ be the canonical
momentum mapping, $\lam\in\gu^{*}$ a regular value in the image of
$\mu$, and $\orb$ the coadjoint orbit passing through $\lam$. Then we
have
\begin{align*}
 (T^{*}U)\spr{\mathcal{O}}G
 &= 
 (T^{*}U)\spr{\lambda}G
 = 
 (T^{*}G\times T^{*}S)\spr{\lam}G
 = 
 (T^{*}G)\spr{\lambda}G\times T^{*}S\\
 &=
 \orb\times T^{*}(U/G)
\end{align*}
as symplectic spaces; since $T^{*}G\spr{\lam}G = \orb$.
The aim of this section is to drop the
freeness assumption. To do so we will take the same approach as Schmah
\cite{Sch04} and use singular commuting reduction.

Now we return to the case where $U=G\times_{H}S$ as introduced
above. On $G\times S$ we will be concerned with two commuting
actions. These are 
\begin{align}
 \lam: G\times G\times S&\longto G\times S, \quad\lam_{g}(k,s)=(gk,s)\\
 \tau: H\times G\times S&\longto G\times S, \quad\tau_{h}(k,s)=(kh^{-1},h.s).
\end{align}
These actions obviously commute. The latter, i.e., $\tau$ is called the
twisted action by $H$ on $G\times S$. We can cotangent lift $\lam$ and
$\tau$ to give Hamiltonian transformations on $T^{*}(G\times S)$ with
momentum mappings $J^{\lam}$ and $J^{\tau}$, respectively. By left
translation we trivialize 
$T^{*}(G\times S) 
 =
(G\times\gu^{*})\times T^{*}S$. 

To facilitate the notation we will denote the cotangent lifted action
of $\lam$, $\tau$ again by $\lam$, $\tau$ respectively. 

\begin{lemma}\label{lem:actions}
Let $(k,\eta;s,p)\in G\times\gu^{*}\times T^{*}S$. Then we have the
following formulas.
\begin{align}
 J^{\lam}(k,\eta;s,p) 
 &= \Ad(k^{-1})^{*}.\eta\label{eq:J-lam}
  =: \Ad^{*}(k).\eta\in\gu^{*},\\
 J^{\tau}(k,\eta;s,p) 
 &= -\eta|\ho+\mu(s,p)\in\ho^{*}\label{eq:J-tau}
\end{align} 
where $\mu$ is the canonical momentum map on $T^{*}S$. Moreover,
the actions $\lam$ and $\tau$ commute, and 
$J^{\lam}$ is $H$-invariant and $J^{\tau}$ is $G$-invariant. 
\end{lemma}

Since the formula of the canonical momentum map on $T^{*}S$ 
with regard to the $H$-action
is the same as that on
$T^{*}U$ 
with regard to the $G$-action
we use the same symbol $\mu$ for  both
these maps. It will be clear from the context whether $\mu$
denotes the $H$- or the $G$-momentum map whence this will not cause any
confusion.  

\begin{proof}
We denote the left action by $G$ on itself by $L$, the right action by
$R$, and the conjugate action by $\conj$. In this notation we then
have $\Ad(k).X=T_{e}\conj_{k}.X$, and
$\conj_{k}=L_{k}\circ R^{k^{-1}}=R^{k^{-1}}\circ L_{k}$. 
It is straightforward to verify that the 
cotangent lifted actions of $L$ and $R$ on $T^{*}G=G\times\gu^{*}$ 
are given by 
\begin{align*}
 T^{*}L_{g}(k,\eta)&=(gk,\eta)=(gk,-\eta\circ\zeta^{R^{-1}}(g))\\
 T^{*}R^{g}(k,\eta)&=(kg,\Ad^*(g^{-1}).\eta)=(kg,\eta\circ\zeta^{L}(g^{-1}))
\end{align*} 
where $\zeta^{L}$ and $\zeta^{R}$ denote the fundamental vector field
mappings associated to $L$ and $R$ respectively.
Using the left trivialization $TG = G\times\gu$ we thus find that 
\begin{align*}
 \vv<J^{\lam}(k,\eta;s,p),X> 
 &= \vv<\eta,\zeta^{L}_{X}(k)> 
 = \vv<\eta,T_kL_{k^{-1}}\dd{t}{}|_0\exp{(tX)}k>\\
 &= \vv<\eta,T_e(L_{k^{-1}}\circ R^k).X>
 = \vv<\Ad^{*}(k).\eta,X>
\end{align*} 
for all $X\in \gu$ which shows the first claim. 
Likewise, it furthermore follows that $\vv<J^{\tau}(k,\eta;s,p),Z> = \vv<-\eta,Z> +
\vv<p,\zeta_{X}(s)>$ for all $Z\in\ho$. 
The invariance of $J^{\lam}$ and $J^{\tau}$ is
immediate from the formulas of the trivialized cotangent lifted
actions. 
\end{proof}

\begin{corollary}
Let $\alpha\in\gu^{*}$ and $\beta\in\ho^{*}$ such that $\alpha$,
$\beta$ is in the image of $J^{\lam}$, $J^{\tau}$ respectively. Then
the following are true.
\begin{enumerate}[\up (1)]
\item The action $\lam$ descends to a Hamiltonian action on the
  Marsden-Weinstein reduced space $T^{*}(G\times
  S)\spr{\beta}H$. Moreover, $J^{\lam}$ factors to a \momap 
  $j_{\lam}: T^{*}(G\times S)\spr{\beta}H\to\gu^{*}$ 
  for this action. 
\item The action $\tau$ descends to a Hamiltonian action on the
  Marsden-Weinstein reduced space $T^{*}(G\times
  S)\spr{\alpha}G$. Moreover, $J^{\tau}$ factors to a \momap 
  $j_{\tau}: T^{*}(G\times S)\spr{\alpha}G\to\ho^{*}$ 
  for this action. 
\item The product action $G\times H\times T^{*}(G\times S)\to
  T^{*}(G\times S)$, $(k,h,u)\mapsto \lam_{k}.\tau_{h}.u$ is
  Hamiltonian with \momap $(J^{\lam},J^{\tau})$. Moreover,
\begin{align*}
 (T^{*}(G\times S)\spr{\alpha}G)\spr{\beta}H
 =
 T^{*}(G\times S)\spr{(\alpha,\beta)}(G\times H)
 =
 (T^{*}(G\times S)\spr{\beta}H)\spr{\alpha}G
\end{align*}
as singular symplectic spaces. 
\end{enumerate}
\end{corollary}

\begin{proof}
Since the actions by $\lam$ and $\tau$ are free the first two
assertions can be deduced from the regular commuting reduction theorem
(\cite{MMOPR03})
with the necessary conditions being verified in the above
lemma. Clearly, the product action by $G\times H$ is well-defined and
Hamiltonian with asserted \momap. However, the product action will not
be free in general. Thus the last point is a consequence of the
singular commuting reduction theorem of Section \ref{s:sg_com_red}.
\end{proof}

We will only be interested in the case where $\beta=0$. 
Moreover,
on $T^*G$ we shall only be concerned with the lifted $\lam$-action. 
Thus the expression
$T^{*}G\spr{\alpha}G$ will throughout stand for
$(J^{\lam})^{-1}(\alpha)/G_{\alpha}$. 

\begin{proposition}\label{prop:bun_pic}
Clearly, $0$ is in the image of $J^{\tau}$. Therefore,
\begin{align*}
 T^{*}U\spr{\alpha}G
 &\cong
  T^{*}(G\times_{H}S)\spr{\alpha}G
 =
  T^{*}(G\times S)\spr{0}H\spr{\alpha}G
 =
  T^*(G\times S)\spr{\alpha}G\spr{0}H\\
 &=
 (T^{*}G\spr{\alpha}G\times T^{*}S)\spr{0}H
 = 
 (\orb\times T^{*}S)\spr{0}H
\end{align*}
as stratified symplectic spaces, and
where $\orb=\Ad^{*}(G).\alpha$. 
\end{proposition}

\begin{proof}
Since the \iso $T^{*}U\xrighto{\simeq}T^{*}(G\times_{H}S)$ comes
from an equivariant \diffeo $U\xrighto{\simeq}G\times_{H}S$ on the base it is
an equivariant \symplecto that intertwines the respective momentum
maps. Now the regular reduction theorem 
for cotangent bundles at zero
momentum says that $T^{*}(G\times_{H}S)$ and $T^{*}(G\times
S)\spr{0}H$ are symplectomorphic. 
(See \cite[Theorem~4.3.3]{AM78} and the remark immediately below \cite[Theorem~4.3.3]{AM78}.)  
Further it is well-known 
(and immediate from Lemma~\ref{lem:actions}(\ref{eq:J-lam}))
that
$T^{*}G\spr{\alpha}G=\orb$. The rest is a direct consequence of
Theorem \ref{thm:sg_comm_red} on singular commuting reduction.  
\end{proof}

From now on we make the assumption that 
\[
 Q=Q_{(H)},
\]
i.e., all isotropy subgroups of points $q\in Q$ are conjugate within
$G$ to $H$. 
Obviously, this assumption imposes a rather strong restriction on the
generality of the subsequent. However, in applications such as in the
Calogero-Moser system of Section~\ref{sec:cms} this is the generic case
in a certain sense. 
See also Remark~\ref{rem:stratQ}.

\begin{theorem}[Bundle picture]\label{thm:bun_pic}
Let $Q=Q_{(H)}$ and
let $\orb\subeq\gu^{*}$ be a coadjoint orbit in the image 
of the \momap 
$\mu: T^{*}Q\to\gu^{*}$. 
Then, locally, we have a singular symplectic fiber bundle 
\[
\xymatrix{
 {\orb\spr{0}H}\ar @{^{(}->}[r] 
 &{T^{*}Q\spr{\mathcal{O}}G} \ar[r]
 &{T^{*}(Q/G)}
}
\]
with typical fiber the singular symplectic space $\orb\spr{0}H$ and smooth base
$T^{*}(Q/G)$. 
\end{theorem}

The fiber bundle in this theorem is singular in the
sense that it is a topological fiber bundle and the transition
functions act by strata preserving transformations on $\orb\spr{0}H$
which are smooth in the sense that they preserve the algebra
$\cinf(\orb\spr{0}H) := W^{\infty}(\orb\cap\ann\ho)^H$
where $W^{\infty}(\orb\cap\ann\ho)$ denotes the Whitney $\cinf$
functions on $\orb\cap\ann\ho\subset\orb$. 
See
\cite{ACG91,Pfl01,Pfl01a} for more on smooth structures on singular
(symplectic) spaces.

\begin{proof}
Consider a tube $U$ of the $G$-action on $Q$. By virtue of the Slice
Theorem (\cite{PT88,DK99,OR04}) 
there thus exists a slice $S$ such that there is a
$G$-equivariant \diffeo
\[
 U\cong G\times_{H}S = G/H\times S.
\]
Indeed, this is true 
since all points of $Q$ are regular by assumption 
whence the slice
representation is trivial. We can lift this diffeomorphism to a
symplectomorphism of cotangent bundles to get 
\[
 T^{*}U\spr{\mathcal{O}}G 
 \cong \orb\spr{0}H\times T^{*}S
\]
as in Proposition \ref{prop:bun_pic} above. Since $T^{*}S$ is a typical
neighborhood in $T^{*}(Q/G)$ the result follows.
\end{proof}

\begin{rem}[On fully singular reduction]\label{rem:stratQ}
For the purpose of this remark assume that the isotropy lattice of the
$G$-action on $Q$ consists of more than one isotropy class. Let $(H)$
be an isotropy class on $Q$, and let $(L)$ be an isotropy class of the
lifted $G$-action on $T^*Q$. 
Let $\ann Q_{(H)}\to Q_{(H)}$ denote the sub-bundle of $(T^*Q)|Q_{(H)}$
consisting of those co-vectors which vanish upon insertion of a vector
tangent to $Q_{(H)}$. Clearly, we have
\[
 (T^*Q)_{(L)}|Q_{(H)}
 =
 (T^*Q_{(H)}\times_{Q_{(H)}}\ann Q_{(H)})_{(L)},
\]
and note that the \momap $\mu: T^*Q\to\gu^*$ vanishes on 
$\ann Q_{(H)}$. 
Therefore,
for an orbit $\orb$ in the image of $\mu$ we have that 
\[
 \mu^{-1}(\orb)|Q_{(H)}
 = 
 \mu_{(H)}^{-1}(\orb)\times_{Q_{(H)}}\ann Q_{(H)}
\]
where $\mu_{(H)}$ 
denotes the \momap of the cotangent lifted
$G$-action on $T^*Q_{(H)}$. 
The $G$-equivariant projection 
$\mu_{(H)}^{-1}(\orb)\times_{Q_{(H)}}\ann Q_{(H)}\to\mu_{(H)}^{-1}(\orb)$
gives rise to a 
mapping 
\[
\tag{F}\label{equ:F}
 \eta_{(L)}:
 (\mu_{(H)}^{-1}(\orb)\times_{Q_{(H)}}\ann Q_{(H)})_{(L)}/G
 \longto\mu_{(H)}^{-1}(\orb)/G
 = T^*(Q_{(H)})\sporb G
\]
the base of which is 
described by Theorem~\ref{thm:WWred} in the presence of a generalized 
connection form on $Q_{(H)}\toto Q_{(H)}/G$.
The map $\eta_{(L)}$ is, in general, neither surjective nor does it have locally constant fiber type.
The fiber over a point $[x]\in\mu_{(H)}^{-1}(\orb)/G$ 
(such that $G_{\tau(x)}=H$)
is of the form 
\[
 \eta_{(L)}^{-1}([x]) 
 = 
 \set{w\in\ann_{\tau(x)}Q_{(H)}: H_w\cap G_x\text{ is conjugate to } L \text{ within } G}/G_x
\] 
where $\tau: T^*Q_{(H)}\to Q_{(H)}$ 
is the cotangent projection.
Note that $gLg^{-1}\subset G_x\subset G_{\tau(x)}$ 
for some $g\in G$
by equivariance of projections.
The image of $\eta_{(L)}$ clearly is a union of orbit type strata. Moreover, using the 
notation of Duistermaat and Kolk~\cite[Definition~2.6.1]{DK99} it is evident
that
\[
 \im\eta_{(L)}\subset(\mu_{(H)}^{-1}(\orb))_{x}^{\lesssim}/G
 = G.(\mu_{(H)}^{-1}(\orb))^L/G
\]
where
$x\in(\mu_{(H)}^{-1}(\orb))_{(L)} = (\mu_{(H)}^{-1}(\orb))_x^{\sim}$, 
whence it follows that 
$(\mu_{(H)}^{-1}(\orb))_{(L)}/G\subset\im\eta_{(L)}$ is open and dense 
since it is the regular stratum of $\im\eta_{(L)}$. 
Let 
$M_0 
 := (T^*(Q_{(H)})\sporb G)_{(L)} 
 := (\mu_{(H)}^{-1}(\orb))_{(L)}/G$ 
(see Theorem~\ref{thm:WWred})
and consider the restriction 
$\eta_0 := \eta_{(L)}|\eta_{(L)}^{-1}(M_0)$ which yields a bundle like object
\[
 \eta_0:
 \eta_{(L)}^{-1}(M_0)\longto M_0,
\]
that is, $\eta_0$ is surjective and the fiber over a point $[x]\in
M_0$ such that $G_x=L$ is $(\ann_{\tau(x)}Q_{(H)})^L$. 
We believe that this object can be shown to constitute a smooth fiber bundle. 
However, the employability of this `bundle'
is quite limited by the fact that
we cannot give a satisfactorily useful description of
$\eta_{(L)}^{-1}(M_0)$. 
Further problems are deciding what a generalized connection form on
$Q\toto Q/G$ should be and determining how the `secondary
strata' 
$(\mu_{(H)}^{-1}(\orb)\times_{Q_{(H)}}\ann Q_{(H)})_{(L)}/G$ fit
together to yield the `primary stratum'
$(\mu^{-1}(\orb)\cap(T^*Q)_{(L)})/G$. The latter problem was solved in
\cite{PRS03} for reduction at trivial orbits $\orb=\set{\textup{point}}$.

If $Q_{(H)}=Q_{\textup{reg}}$ is the regular stratum which is open dense in
$Q$ then $\ann Q_{(H)}$ is trivial. In this (generic) case
Theorems~\ref{thm:bun_pic} and \ref{thm:WWred} thus provide a full answer to the
reduction problem. In the more general situation these results clearly
provide only a partial answer to the reduction problem. However, it is
expected that any solution to this problem will rely on these single
orbit type results. 
\end{rem}

\section{Gauged cotangent bundle reduction}\label{sec:gauged_red}

Continue
to assume that we are in the situation of Section~\ref{sec:bun-pic}.
In particular, we suppose that $Q=Q_{(H)}$ is of single orbit type.
However, as an additional input datum we assume from now on a
generalized
principal bundle connection form $A\in\Om^1(Q;\gu)$ on $Q\toto Q/G$
given. The term generalized is to be understood in the context of
Alekseevsky and Michor~\cite[Section~3.1]{AM95}. This means that $A: TQ\to\gu$ is
$G$-equivariant and that $\zeta = \zeta\circ A\circ\zeta$. 
In particular, the connection form $A$ induces a right inverse to the
projection $\gu\toto\gu/\gu_q$ depending smoothly on $q\in Q$. 

According to
\cite[Section 4.6]{AM95} the curvature form associated to $A$ is
defined by 
\[
 \curv^{A}
 :=
 dA-\by{1}{2}[A,A]^{\wedge}
\]
where 
\[
 [\var,\psi]^{\wedge}(v_1,\dots,v_{l+k})
 :=
 \by{1}{k!l!}
  \sum_{\sigma}\sign\sigma
  [\var(v_{\sigma1},\dots,v_{\sigma l}),\psi(v_{\sigma(l+1)},\dots,v_{\sigma(l+k)})]
\]
is the graded Lie bracket on 
$
 \Om(Q;\gu) 
 := \bigoplus_{k=0}^{\infty}\Gamma(\Lam^{k}T^{*}Q\otimes\gu)
$, 
and $\var\in\Om^{l}(Q;\gu)$
and $\psi\in\Om^{k}(Q;\gu)$. The sign in our definition of $\curv^A$
differs from that in \cite{AM95} because we are concerned with left
$G$-actions as opposed to right actions.

Since the $G$-action on $Q$ is of single orbit type the orbit
space $Q/G$ is a smooth manifold, and the projection $\pi: Q\to Q/G$
is a fiber bundle with typical fiber $G/H$. 
However, the isotropy lattice of the lifted action by $G$ on
$T^{*}Q$ is, in general (for $H\neq\set{e}$), non-trivial whence
the quotient space $(T^{*}Q)/G$
is a stratified space. Its strata are of the form
$(T^{*}Q)_{(L)}/G$ where $(L)$ is in the isotropy lattice of $T^{*}Q$. 

The vertical sub-bundle of $TQ$ with respect to $\pi: Q\to Q/G$ is
$\ver := \ker T\pi$. Via the connection $A$ we can also define the 
horizontal sub-bundle $\hor := \ker A$. We define the dual
horizontal sub-bundle of $T^{*}Q$ as the sub-bundle $\hor^{*}$
consisting of those co-vectors that vanish on all vertical
vectors. Likewise, we define the dual
vertical sub-bundle of $T^{*}Q$ as the sub-bundle $\ver^{*}$
consisting of those co-vectors that vanish on all horizontal
vectors. As usual, the connection $A$ provides a trivialization of the
vertical sub-bundle, i.e., $\ver \cong_{A} \bsc_{q\in Q}\gu/\gu_q$. 
In particular, $\bsc_{q\in Q}\gu/\gu_q$ and $\bsc_{q\in Q}\ann\gu_q$
are smooth vector bundles.

\subsection{Mechanical connection}\label{sub:mech-con}
If $(Q,\vv<.,.>)$ is a Riemannian manifold and $G$ acts on $Q$ by isometries 
there is a certain connection which is particularly well adapted to
mechanical systems on $Q$. This is the so-called mechanical connection
which is defined as follows. 
For $X,Y\in\gu$ and $q\in Q$ we define
$\ine_{q}(X,Y) := \vv<\zeta_{X}(q),\zeta_{Y}(q)>$ and call this  
the \caps{locked inertia tensor}.  
This defines a
non-degenerate pairing on $\gu/\gu_q$ whence
it provides an identification $\check{\ine_{q}}:
\gu/\gu_{q}\to(\gu/\gu_{q})^{*} = \ann\gu_{q}$. We use this \iso
to define a one-form $\tilde{A}$ on $Q$ with values in the bundle 
$\bsc_{q\in Q}\gu/\gu_{q}$ by the following diagram.
\[
\xymatrix{
 {T_{q}^{*}Q} \ar[r]^-{\mu_{q}} 
              &
 {\ann\gu_{q}}
 \ar[d]^{(\check{\mathbb{I}_{q}})^{-1}}
             \\
 {T_{q}Q}     \ar[u]^{\simeq}
              \ar @{-->}[r]^{\tilde{A}_{q}}& 
 {\gu/\gu_{q}}
}
\]
Notice that $\im\mu_q = \ann\gu_q$ by reason of dimension. 
Thus we have a trivialization 
$\ver\cong_{\tilde{A}}\bsc_{q\in Q}\gu/\gu_q$ 
of the vertical sub-bundle. However, to obtain a
generalized connection form on $Q\toto Q/G$ from this trivialization
we have to assume one additional object: namely, let $r_q: \gu/\gu_q$
be a right inverse to the projection $\gu\toto \gu/\gu_q$ depending
smoothly on $q\in Q$. 
The (generalized) \caps{mechanical connection} $A: TQ\to\gu$ 
on $Q\toto
Q/G$ is thus defined as the composition 
$A_q = r_q\circ\tilde{A}_q: 
 T_qQ\to\gu/\gu_q\to\gu$.

One obvious way to obtain such a right inverse
is to choose a $G$-invariant non-degenerate bilinear
form on $\gu$ such that $\gu/\gu_q\cong\gu_q^{\bot}\hookto\gu$ with respect to
this form. In many examples such a non-degenerate form is given
canonically.

The mechanical connection was first defined by Smale \cite{Sma70} 
for the case of Abelian group actions. See also
Marsden, Montgomery, and Ratiu \cite[Section 2]{MMR90}. 

To verify that the mechanical connection is indeed a generalized
principal connection form one checks that
$A: TQ\to\gu$
is equivariant
and 
$\zeta(q)(A_{q}(\zeta_{X}(q)))=\zeta_X(q)$ for all $X\in\gu$.

\subsection{Weinstein realization of $T^*Q$}
Let $A$ continue to denote a  generalized
principal connection form on $Q$, and let $\vv<.,.>$ denote the dual pairing.
We define a point-wise dual 
$A_{q}^{*}: 
 \ann\gu_{q}\to
 \ver_{q}^{*}\subeq T_{q}^{*}Q$ 
by the formula 
$
 \vv<A_{q}^{*}(\lam),v> 
 = 
 \vv<\lam,A_{q}(v)>
$
where $\lam\in\ann\gu_{q}$ and $v\in T_{q}Q$. Notice that 
$
 A_{q}^{*}(\mu_{q}(p)) = p
$
for all $p\in\ver_{q}^{*}$
and 
$
 \mu_{q}(A_{q}^{*}(\lam)) = \lam
$
for all $\lam\in\ann\gu_{q}$ since $A$ is a connection form. 

Let $\pi: Q\toto Q/G$ and $\tau_Q: TQ\to Q$ denote the
projections. From the connection form $A$ we obtain the horizontal
lift mapping which we denote by 
\[
 C := ((\tau_Q,T\pi)|\hor)^{-1}:
 Q\times_{Q/G}T(Q/G)\to \hor\hookto TQ. 
\]
Its fiber restriction shall be
denoted by $C_q: \set{q}\times T_{\pi(q)}(Q/G)\to\hor_q\hookto T_qQ$. 

Using the horizontal lift
$C$
on the one hand and the connection
$A$ on the other hand 
we obtain a $G$-equivariant \iso 
\[
 TQ = \hor\oplus\ver 
 \longto
  (Q\times_{Q/G}T(Q/G))\times_{Q}\bsc_{q\in Q}\gu/\gu_{q}
\]
of bundles over $Q$. 
There is  a dual 
version to this \iso, and we will abbreviate 
\[
 \WW :=
 (Q\times_{Q/G}T^{*}(Q/G))\times_{Q}\bsc_{q\in Q}\ann{\gu_{q}}
 \cong
 \hor^{*}\oplus\ver^{*}.
\]
The explicit form of the \iso $\WW\cong T^*Q$ is stated in
Proposition~\ref{prop:WW} below.
To set up some notation for the upcoming proposition, 
and clarify the picture consider the following
stacking of pull-back diagrams.
\[
\xymatrix{
 {\WW} \ar[r]^-{\rho^{*}\widetilde{\tau}=\widetilde{\widetilde{\tau}}}
       \ar[d]_{\widetilde{\tau}^{*}\rho=\widetilde{\rho}} &
 {\bsc_{q}\ann\gu_{q}}
       \ar[d]^{\rho} \\ 
 {Q\times_{Q/G}T^{*}(Q/G)}
       \ar[r]^-{\pi^{*}\tau=\widetilde{\tau}}
       \ar[d]_{\tau^{*}\pi=\widetilde{\pi}} &
 Q 
       \ar[d]^{\pi}  \\
 T^{*}(Q/G)  
       \ar[r]^-{\tau} &
 Q/G       
}
\]
The upper stars in this diagram are, of course, not pull-back
stars. By slight abuse of notation we shall denote elements
$(q;\pi(q),\eta;q,\lam)\in\WW$ simply by $(q,\eta,\lam)$.  
Further, let $\tau_{\mathcal{W}}: T\WW \to \WW$ denote the tangent
projection.

\begin{proposition}[Symplectic structure on \WW]\label{prop:WW}
The chosen connection form $A$ induces an \iso
\begin{align*}
 \psi=\psi(A):
 (&Q\times_{Q/G}T^{*}(Q/G))\times_{Q}\bsc_{q\in Q}\ann\gu_{q} = \WW
  \longto T^{*}Q, \\
 (&q,\eta,\lam)
    \longmapsto
   (q,(T_q\pi)^*\eta+A_q^{*}(\lam))
\end{align*}
of bundles over $Q$,
and the following are true.
\begin{enumerate}[\up (1)]
\item
There is an induced $G$-action on $(\WW,\sigma=\psi^*\Om)$ by 
symplectomorphisms.
Here $\Om=-d\theta$ is the canonical symplectic form on $T^*Q$.
Moreover, this action is Hamiltonian with momentum
map
\[
 \mu_{A} = \mu\circ\psi: \WW\longto\gu^{*}, (q,\eta,\lam)\longmapsto\lam,
\]
where $\mu$ is the
momentum map $T^{*}Q\to\gu^{*}$, 
and $\psi$ is equivariant.
\item
The induced symplectic form $\sigma$ on the connection
dependent realization $\WW$ of $T^{*}Q$ 
is given by the formula
\[
 \sigma =
 (\wt{\pi}\circ\wt{\rho})^{*}\Om^{Q/G} - dB.
\]
Here
$\Om^{Q/G}=-d\theta^{Q/G}$ is the canonical symplectic form on $T^{*}(Q/G)$, 
and 
$B\in\Om^{1}(\WW)$ 
is given by 
\[
 B
 = 
 \vv<\tau_{\mathcal{W}}^*\mu_A,(\wt{\tau}\circ\wt{\rho})^*A>.
\]
Moreover,
\begin{align*}
 dB 
 &=
  \vv<\tau_{\mathcal{W}}^*d\mu_A\stackrel{\wedge}{,}(\wt{\tau}\circ\wt{\rho})^*A>
  + \vv<\tau_{\mathcal{W}}^*\mu_A , (\wt{\tau}\circ\wt{\rho})^*\curv^A>\\
 &\phantom{=}
  + \vv<\tau_{\mathcal{W}}^*\mu_A , \by{1}{2}[(\wt{\tau}\circ\wt{\rho})^*A,(\wt{\tau}\circ\wt{\rho})^*A]^{\wedge}>.
\end{align*} 
Here $\vv<.\stackrel{\wedge}{,}.>$ denotes the exterior multiplication
of a $\gu^*$-valued form with a $\gu$-valued form. 
In local coordinates where we may use a splitting of tangent vectors
$\xi_1,\xi_2\in T_{(q,\eta,\lam)}\WW$ as 
$\xi_i 
 = (q'_i,\eta'_i,\lam'_i)
$
and 
$q'_i = v_{i}^{\textup{hor}}+\zeta_{Z_{i}}(q)$
for $i=1,2$
this means that 
\[
 dB_{(q,\eta,\lam)}(\xi_1,\xi_2) \\
  = 
 \vv<\lam'_1,Z_2> - \vv<\lam'_{2},Z_1> 
 + 
 \vv<\lam,\curv^{A}_{q}(q'_1,q'_2)>
  + \vv<\lam,[Z_{1},Z_{2}]>.
\]
\end{enumerate}
\end{proposition}

\begin{proof}
Obviously, $\psi$ is a bundle map. Its inverse is given by 
the bundle map
$(q,p)\mapsto(q,C_q^*(p),\mu_q(p))$ where $C_q: \set{q}\times
T_{\pi(q)}(Q/G)\to T_qQ$ is the horizontal lift mapping and $\mu_q :=
\mu|T_q^*Q$ is the fiber restriction of the \momap. 

Assertion (1) is clear from the construction.

To see assertion (2) we work locally in $Q$. That is, let $U$ be a
trivializing patch for $\WW\to Q$ as well as for $T^*Q\to Q$ and
consider $\psi|U: \WW|U\to T^*Q|U$. Let $(q,\eta,\lam)\in\WW|U$ and
$\xi\in T_{(q,\eta,\lam)}\WW$. By locality we may split $\xi$ as $\xi
= (q',\eta',\lam')$ where $q'\in T_qQ$, $\eta'\in T_{\eta}(T^*(Q/G))$,
and $\lam'\in T_{\lam}(\bsc_{q\in Q}\ann\gu_q)$. 
However, in order to make the notation not too cumbersome we do not invent new
symbols (like $\xi^U$ or $A^U$) for the local versions of global
objects (like $\xi$ or $A$). 
To find the desired
formula for $\sigma = \psi^*\Om = -d\psi^*\theta$, note that
\begin{align*}
 (\psi^*\theta)_{(q,\eta,\lam)}(\xi) 
 &=
 \theta_{(q,(T_q\pi)^*\eta+A_q^*(\lam))}(T_{(q,\eta,\lam)}\psi.\xi)\\
 &=
 \theta_{(q,(T_q\pi)^*\eta+A_q^*(\lam))}(T_{(q,\eta,\lam)}\psi^1.\xi,T_{(q,\eta,\lam)}\psi^2.\xi)\\
 &=
 \vv<(T_q\pi)^*\eta+A_q^*(\lam),T_{(q,(T_q\pi)^*\eta+A_q^*(\lam))}\tau_Q.T_{(q,\eta,\lam)}\psi^1.\xi>\\
 &=
 \vv<(T_q\pi)^*\eta,q'> + \vv<A_q^*(\lam),q'>\\
 &=
 \vv<\eta,T_q\pi.q'> + \vv<\lam,A_q(q')>\\
 &=
 ((\wt{\pi}\circ\wt{\rho})^*\theta^{Q/G})_{(q,\eta,\lam)}(\xi)
 +
 (\vv<\tau_{\mathcal{W}}^*\mu_A,(\wt{\tau}\circ\wt{\rho})^*A>)_{(q,\eta,\lam)}(\xi)\\
 &=
 (((\wt{\pi}\circ\wt{\rho})^*\theta^{Q/G})+B)_{(q,\eta,\lam)}(\xi)
\end{align*}
where 
$\psi|U = (\psi^1,\psi^2): \WW|U\to U\times V = T^*Q|U$ 
and $V$ is the standard fiber of 
$\tau_Q: T^*Q\to Q$. 
Since the first
and the last expressions in this computation are global objects it is
true that
\begin{align*}
 \sigma 
 &=
 (\wt{\pi}\circ\wt{\rho})^*\Om^{Q/G}
 - 
 dB\\
 &=
 (\wt{\pi}\circ\wt{\rho})^*\Om^{Q/G}
  - 
    \vv<\tau_{\mathcal{W}}^*d\mu_A\stackrel{\wedge}{,}(\wt{\tau}\circ\wt{\rho})^*A>\\
 &\phantom{=i} 
  - \vv<\tau_{\mathcal{W}}^*\mu_A , (\wt{\tau}\circ\wt{\rho})^*\curv^A>
  + \vv<\tau_{\mathcal{W}}^*\mu_A , \by{1}{2}[(\wt{\tau}\circ\wt{\rho})^*A,(\wt{\tau}\circ\wt{\rho})^*A]^{\wedge}>.
\end{align*}
Indeed, this is because  
$\curv^A = dA-\by{1}{2}[A,A]^{\wedge}$.
\end{proof}

The $G$-action on $\WW$ is, of course, given by $g.(q,\eta,\lam) = (g.q,\eta,\Ad^*(g).\lam)$. 
Similarly there is an induced $G$-action on 
$\bsc_{q\in Q}\gu/\gu_q \cong \ver$ 
which is given by $g.(q,X+\gu_q) = (g.q,\Ad(g).X+\gu_{g.q})$. 

The notion of a stratified map which appears in the following theorem
is defined in Section~\ref{sec:prel}.

\begin{theorem}[Weinstein space]
There are stratified isomorphisms of stratified bundles over $Q/G$:
\begin{align*}
 \alpha=\alpha(A): 
 \bsc_{(L)}(TQ)_{(L)}/G &\longto 
 T(Q/G)\times_{Q/G}\bsc_{(L)}(\bsc_{q\in Q}\gu/\gu_{q})_{(L)}/G,\\
 [(q,v)] &\longmapsto(T\pi(q,v),[(q,A_{q}(v))])
\end{align*}
where $(L)$ runs through the isotropy lattice of $TQ$. 
The dual \iso is given by
\begin{align*}
 \beta=(\alpha^{-1})^{*}: 
 (T^{*}Q)/G &\longto T^{*}(Q/G)\times_{Q/G}(\bsc_{q\in
                                Q}\ann\gu_{q})/G =: W,\\
 [(q,p)] &\longmapsto(C^{*}_q(p),[(q,\mu(q,p))])
\end{align*}
where the stratification was suppressed. 
Here 
\[
\smash{
 C^{*}_q: 
 T^{*}_qQ\overset{\iota_q^*}{\to}\hor_q^{*}\to T_{\pi(q)}^{*}(Q/G)
}
\] 
is 
the point wise dual to the horizontal lift mapping 
\[
\smash{
C:
Q\times_{Q/G}T(Q/G)\longto\hor\overset{\iota_q}{\hookto} TQ}, 
 \text{ }
([q],v;q)\longmapsto C_{q}(v). 
\]
Moreover, $\beta$ is an \iso of Poisson spaces as follows: we can
naturally identify 
\begin{align*}
 \WW/G \overset{=}{\longto} W, \;
 [(q;[q],\eta;q,\lam)] \longmapsto ([q],\eta;[(q,\lam)])
\end{align*}
thus obtaining a quotient Poisson bracket on
$\cinf(W)=\cinf(\WW)^{G}$ as the quotient Poisson bracket.
\end{theorem} 

In the case that $G$ acts on $Q$ freely the first assertion of the
above theorem  can also be
found in Cendra, Holm, Marsden, Ratiu \cite{CHMR98}. Following Ortega
and Ratiu \cite[Section 6.6.12]{OR04} the above constructed
interpretation $W$ of $(T^{*}Q)/G$ is called \caps{Weinstein space}
referring to Weinstein \cite{Wei78} where this universal
construction first appeared. 

In fact, the original construction of \cite{Wei78} was the following:
Let $Q$ be a left free and proper $G$-space such that $Q\toto Q/G$ is
endowed with a principal bundle connection form $A$ and let $F$ be a
right Hamiltonian $G$-space with equivariant \momap $\Phi: F\to\gu^*$. 
To make $F$ into a left Hamiltonian $G$-space we use the inversion in
the group. The \momap with respect to the thus obtained $G$-action is
given by $-\Phi$.  
Under these assumptions \cite{Wei78} proves
that the smooth symplectic quotient 
$
 (T^*Q\times F)\spr{0}G 
 =
 (\mu-\Phi)^{-1}(0)/G
 \cong_A
 T^*(Q/G)\times_{Q/G}(Q\times_G F)
$  
is symplectomorphic to the Sternberg space
$(Q\times_{Q/G}T^*(Q/G))\times_G F$ of \cite{Ste77}. Taking $F$ to be
the Hamiltonian $G$-space 
$\orb$ acted upon by $\Ad^*(G)$, and employing the shifting trick
$T^*Q\sporb G \cong (T^*Q\times\orb)\spr{0}G$, this construction
yields the realization 
$T^*Q\sporb G 
\cong_A
T^*(Q/G)\times_{Q/G}(Q\times_G \orb)$. In particular, one thus obtains
a fiber bundle $\orb\hookto T^*Q\sporb G \to T^*(Q/G)$.    

It thus makes sense to refer to the realization $W$ of $(T^*Q)/G$ which is
constructed along similar lines as a Weinstein space. The induced Poisson
structure on $W$ is explicitly described in \cite{HR05}.

\begin{proof}
As already noted above $(TQ)/G$ is a stratified space. Since the base
$Q$ is stratified as consisting only of a single stratum, the
equivariant foot
point projection map $\tau: TQ\to Q$ is trivially a stratified
map. Thus, we really get a stratified bundle $(TQ)/G\to Q/G$. In the
same spirit $(\bsc_{q\in Q}\gu/\gu_{q})/G$ is stratified into
orbit types, and the projection onto $Q/G$ is a stratified bundle
map. According to Davis \cite{Dav78} pullbacks are well defined in
the category of stratified spaces and stratified maps
and thus it makes sense to define
$T(Q/G)\times_{Q/G}(\bsc_{q\in Q}\gu/\gu_{q})/G$. 

The map $\alpha$ is well defined: indeed, for $(q,v)\in TQ$ and 
$k\in G$ 
we have 
$T\pi(k.q,k.v) 
 = (\pi(k.q),T_{k.q}\pi(T_{q}l_{k}(v))) 
 = (\pi(q),T_{q}(\pi\circ l_{k})(v)) 
 = T\pi(q,v)$, 
and
$[(k.q,A_{k.q}(k.v))] = [(q,A_q(v))]$ 
by equivariance of $A$.
It is clearly continuous
as a composition of continuous maps. 

We claim that $\alpha$ maps strata onto strata, and moreover we
have the formula
\[
 \alpha((TQ)_{(L)}/G)
 =
 T(Q/G)\times_{Q/G}(\bsc_{q\in Q}\gu/\gu_{q})_{(L)}/G.
\]
Indeed, this follows immediately since $\alpha$ lifts to a smooth equivariant 
\iso
$\wt{\alpha}: TQ\to(Q\times_{Q/G}T(Q/G))\times_Q\bsc_{q\in Q}\gu/\gu_q$, 
$(q,v)\mapsto(q;T\pi(q,v);q,A_q(v))$
of vector bundles over $Q$,
and we have clearly that 
$\wt{\alpha}((TQ)_{(L)}) =  
(Q\times_{Q/G}T(Q/G))\times_Q(\bsc_{q\in Q}\gu/\gu_q)_{(L)}$. 
The restriction of $\alpha$ to any stratum is smooth
as a composition of smooth maps. 

Since $\zeta(q)(A_{q}(\zeta_{X}(q))) = \zeta_X(q)$ for
$X\in\gu$ we can write down an inverse 
\[
 \alpha^{-1}:
 ([q],v;[(q,X)])
 \to
 [(q,C_{q}(v)+\zeta_{X}(q))]
\]
and again it is an easy matter to notice that this map is well
defined, continuous, and smooth on each stratum.
In fact, we have used here that the connection $A$ by definition
provides a right inverse to $\gu\toto\gu/\gu_q$ whence by slight
abuse of notation we may consider elements $X\in\gu/\gu_q$ as elements
in $\gu$.

It makes sense to define the dual $\beta$ of the inverse map
$\alpha^{-1}$ in a point wise manner, and it only remains to compute
this map.
\begin{align*} 
 \vv<\beta[(q,p)],([q],v;[(q,X)])> 
 &= 
 \vv<[(q,p)],[(q,C_{q}(v)+\zeta_{X}(q))]> \\
 &= 
 \vv<p,C_{q}(v)> + \vv<p,\zeta_{X}(q)> \\
 &=
 \vv<C^{*}_q(p),v> + \vv<\mu(q,p),X> \\
 &=
 \vv<(C^{*}_q(p),[(q,\mu(q,p))]),([q],v;[(q,X)])>
\end{align*}
where we used the $G$-invariance of the dual pairing over $Q$. 

Finally, $\beta$ is an \iso of Poisson spaces: 
note first that the identifying map 
$\WW/G\to W$, 
$[(q;[q],\eta;q,\lam)]\mapsto ([q],\eta;[(q,\lam)])$ 
is well-defined because $G_{q}$ acts
trivially on $\hor_{q}^{*} \cong T^{*}_{[q]}(Q/G)\ni\eta$ which in turn is due
to the fact that all points of $Q$ are regular. The quotient Poisson
bracket is well-defined since $\cinf(\WW)^{G}\subeq\cinf(\WW)$ is a
Poisson sub-algebra. The statement now follows because the diagram
\[
\xymatrix{
 T^{*}Q \ar[r]^-{\psi^{-1}}\ar @{->>}[d] & {\WW}\ar @{->>}[dr] \\
 (T^{*}Q)/G\ar[r]^-{\beta} & W\ar @{=}[r] & {\WW}/G
}
\]
is commutative, and composition of top and down-right arrow is
Poisson and the left vertical arrow is surjective. 
\end{proof}

\subsection{The reduced phase space}
The following lemmas are key to the subsequent. They
guarantee, in particular, 
that every non-empty pre-image of a coadjoint orbit
$\orb\subset\gu^*$ under $\mu$ fibrates surjectively over $\hor^*$. 
In Theorem~\ref{thm:WWred} we use this to show that every non-empty 
symplectic stratum of the 
reduced phase space fibrates over $T^*(Q/G)$.
The meta principle motivating these results is that $\mu$ 
(which is defined by means of the universal connection $\theta$)
can be thought of 
as a kind of
universal connection form on $Q\toto Q/G$.

\begin{lemma}
Let $\orb\subeq\gu^{*}$ be a coadjoint orbit, $\mu: T^{*}Q\to\gu^{*}$
the canonical momentum mapping, and $\mu_{q}:=\mu|(T_{q}^{*}Q)$. Then
either $\mu_{q}^{-1}(\orb)=\emptyset$ for all $q\in Q$ or
$\mu_{q}^{-1}(\orb)\neq\emptyset$ for all $q\in Q$. In the latter case
we have
\[
 \mu_{q}^{-1}(\orb)
 =
 \ann_{q}(T_{q}(G.q))\times\set{A_{q}^{*}(\lam):
                                \lam\in\ann\gu_{q}\cap\orb}
\]
which is an equality of topological spaces and
where $A_{q}^{*}: \ann\gu_q\to\ver_q^*$ 
is the adjoint of $A_{q}: T_{q}Q\to\gu/\gu_{q}$. 
\end{lemma}

\begin{proof}
Assume firstly that $q_1,q_2\in Q$ lie in the same $G$-orbit. Then it
is obviously true that $\mu_{q_1}^{-1}(\orb)$ is empty if and only if
$\mu_{q_2}^{-1}(\orb)$ is empty. Thus for the purpose of this proof we
can assume that $G_{q_1}=G_{q_2}=H$:
all isotropy subgroups are conjugate to each other, and $q_{2}$ can be
moved around in its orbit without loss of generality.
Now assume that 
$\mu_{q_1}^{-1}(\orb)$ is non-empty, i.e., there is 
$\lam = \mu_{q_1}(p_1)\in\ann\ho\cap\orb$. Using the
connection $A$ we may then define 
$p_2 := A_{q_2}^*(\lam)\in\ver_{q_2}^*\cap\mu_{q_2}^{-1}(\orb)$
whence $\mu_{q_2}^{-1}(\orb)$ is non-empty as well. 
In fact, this construction also proves the last claim of the lemma.
\end{proof}

\begin{lemma}\label{lm:WWorb}
Let $\orb$ be a coadjoint orbit in the image of the \momap $\mu_{A}:
\WW\to\gu^{*}$. Further, let $(L)$ be in the isotropy lattice of the
$G$-action on $\WW$ such that
$\mu_{A}^{-1}(\orb)\cap\WW_{(L)}\neq\emptyset$. Then 
\[
 \WW_{(L)}
 = 
 (Q\times_{Q/G}T^{*}(Q/G))\times_{Q}(\bsc_{q\in Q}\ann\gu_{q})_{(L)}
\]
and 
\[
 \WW_{(L)}\cap\mu_{A}^{-1}(\orb)
 =
 (Q\times_{Q/G}T^{*}(Q/G))\times_{Q}(\bsc_{q\in Q}\ann\gu_{q}\cap\orb)_{(L)}
\]
are smooth manifolds. 
Moreover, 
\[
\xymatrix{
 {\orb_{(L_0)_H^G}\cap\ann\ho}
 \ar @{^{(}->}[r]&
 {(\bsc_{q\in Q}\orb\cap\ann\gu_{q})_{(L)}}
 \ar[r]
 &Q
}
\]
is a smooth fiber bundle where $L_{0}$ is a subgroup of $H$ such that
$L_{0}$ is conjugate to $L$ within $G$.
\end{lemma}

Notice that we do not assume $\orb\cap\ann\ho$ to be smooth.
The notation $\orb_{(L_0)_H^G}$ is explained in Section~\ref{sec:prel}(\ref{def:weirdo-type}).
The $G$-action on $\bsc_{q\in Q}\orb\cap\ann\gu_q$ is induced from the $G$-action on $\WW$
from Proposition~\ref{prop:WW}, and is given by $g.(q,\lam) = (g.q,\Ad^*(g).\lam)$. 

\begin{proof}
The statement about $\WW_{(L)}$ is clear. Thus also the description
of $\WW_{(L)}\cap\mu_{A}^{-1}(\orb)$ follows from the previous lemma
together with Theorem~\ref{thm:sing_spr}.

Concerning the second assertion let $q_{0}\in Q$ with $G_{q_{0}}=H$. Then
\[
 (q_{0},\lam)\in(\bsc_{q\in Q}\orb\cap\ann\gu_{q})_{(L)}
\]
if and only if
\[
 \lam\in\orb\cap\ann\ho
 \textup{ and }
 H\cap G_{\lambda}=H_{\lambda}=L_{0}\textup{ is conjugate to } L
 \textup{ in }G
\]
which is true if and only if
\[
 \lam\in(\orb\cap\ann\ho)_{(L_0)_H^G}
\] 
where $L_{0}$ is a subgroup of $H$ conjugate to $L$ within $G$,
and we view $\orb\cap\ann\ho\subset\orb$ as an $H$-space by
virtue of the restricted $\Ad^*(H)$-action.  
By Lemma~\ref{lem:weirdo-type} the space
$(\orb\cap\ann\ho)_{(L_0)_H^G}$ is a smooth manifold.

To see smooth local triviality we proceed as follows. Let again
$q_{0}\in Q$ with $G_{q_{0}}=H$, and let $S$ be a slice at $q_{0}$ and
$U$ a tube around $G.q_{0}$. That is, $G/H\times S\cong U$,
$(kH,s)\mapsto k.s$ 
as proper $G$-spaces 
by virtue of the Slice Theorem (\cite{PT88,DK99}). 
Then we consider the smooth trivializing map
\begin{align*} 
 S\times G\times_{H}(\orb\cap\ann\ho)_{(L_0)_H^G}
 &\longto
 (\bsc_{q\in Q}\orb\cap\ann\gu_{q})_{(L)}|U,\\
 (s,[(k,\lam_0)]_H)
 &\longmapsto
 (k.s,\Ad^{*}(k).\lam_0)
\end{align*}
which is well defined since 
$[(k,\lam_0)]_H\in G\times_H(\orb\cap\ann\ho)_{(L_0)_H^G}$ implies that
$H_{\lam_0} = gL_0g^{-1}\in (L)$ for some $g\in G$ and this yields
$G_{(k.s,\textup{Ad}^*(k).\lam_0)} = kG_{(s,\lam_0)}k^{-1} 
 = kH_{\lam_0}k^{-1} = kgL_0g^{-1}k^{-1}$.  
Hereby we use, firstly, that 
the diagonal $H$-action cancels
out, i.e., $\Ad^*(kh^{-1}).\Ad^*(h).\lam_0 = \Ad^*(k).\lam_0$ for all
$h\in H$;
secondly, we use that $\gu_s = \gu_{q_0} = \ho$ for all $s\in S$ since $S$ is a 
slice at $q_0$ -- see, e.g., \cite[Corollary~5.1.13(2)]{PT88}.
Clearly, this map is smooth with smooth inverse
$
 (q,\lam)=(k.s,\Ad^{*}(k).\lam_{0})
 \longmapsto
 (s,[(k,\lam_{0})]_H).
$
Therefore, this construction provides smooth bundle charts of the total space 
$(\bsc_{q\in Q}\orb\cap\ann\gu_{q})_{(L)}$.
\end{proof}

\begin{lemma}\label{lem:weirdo-type}
Let $G$, $\orb$, $H$ and $L_0$ be as in Lemma~\ref{lm:WWorb}. Then
$(\orb\cap\ann\ho)_{(L_0)_H^G}$ is a smooth (possibly disconnected) pre-symplectic
manifold. 
\end{lemma}

This fact is an instance of the general theory of
singular commuting reduction of \cite{MMOPR03}. However,
the argument is rather explicit.

\begin{proof}
As in Lemma~\ref{lem:actions} consider the $G$-action on
$T^*G=G\times\gu^*$ given by $g.(k,\lam) = (gk,\lam)$ and the
$H$-action given by $h.(k,\lam) = (kh^{-1},\Ad^*(h).\lam)$. These
actions are Hamiltonian.
By Lemma~\ref{lem:actions}(\ref{eq:J-tau}) 
symplectic reduction of
$G\times\gu^*$ with respect to $H$ at $0\in\ho^*$ yields
$(G\times\gu^*)\spr{0}H = G\times_{H}\ann\ho$. 
Let $\pi_H: G\times\ann\ho\toto G\times_H\ann\ho$ denote the orbit projection.
The corollary of
Lemma~\ref{lem:actions} (or straightforward computation) implies that there
is an induced $G$-action on $G\times_H\ann\ho$ 
given by $g.[(k,\lam)]_H = [(gk,\lam)]_H$
with \momap $j$ given
by $j[(k,\lam)]_H = \Ad^*(k).\lam$. 
(The formula for $j$ follows from
Lemma~\ref{lem:actions}(\ref{eq:J-lam}).)
Now the theory of singular symplectic reduction 
(see Theorem~\ref{thm:sing_spr}) 
implies that 
\begin{multline*}
(j^{-1}(\orb))_{(L_0)^G} 
 = j^{-1}(\orb)\cap(G\times_H\ann\ho)_{(L_0)^G} = \\
   \set{[(k,\lam)]_H\in G\times_H(\orb\cap\ann\ho): 
          G_{[(k,\lam)]_H} = kH_{\lam}k^{-1}
          \textup{ is conjugate to } L_0 \textup{ in } G}
\end{multline*}
is a sub-manifold of
$G\times_H\ann\ho$. Therefore, the pre-image
$\pi_H^{-1}((j^{-1}(\orb))_{(L_0)^G})
 = G\times(\orb\cap\ann\ho)_{(L_0)_H^G}$ is a
sub-manifold of $G\times\ann\ho$, and  
$(\orb\cap\ann\ho)_{(L_0)_H^G}$ is a sub-manifold of $\ann\ho$. 

Finally, the manifold in question is pre-symplectic since
$(\orb\cap\ann\ho)_{(L_0)_H^G}/H 
 \cong (G\times_H(\orb\cap\ann\ho)_{(L_0)^G_H})/G
 = (j^{-1}(\orb))_{(L_0)^G}/G$ is a
symplectic manifold. 
\end{proof}

The singular reduction diagram of Ortega and Ratiu~\cite[Theorem 8.4.4]{OR04}
adjoined to the universal reduction procedure of Arms, Cushman, and
Gotay~\cite{ACG91} (see also \cite[Section 10.3.2]{OR04}) applied to
the Weinstein space has the following form. 
\[
\xymatrix{
 {\mu_{A}^{-1}(\orb)}
   \ar @{<-^{)}}[r]
   \ar @{->>}[d]&
 {\mu_{A}^{-1}(\lam)}
   \ar @{^{(}->}[r]
   \ar @{->>}[d]&
 {\WW}
   \ar @{->>}[d]\\
 {\mu_{A}^{-1}(\orb)/G}
   \ar @{<-}[r]^-{\simeq}&   
 {\mu_{A}^{-1}(\lam)}/G
   \ar @{^{(}->}[r]&
 {\WW/G}
   \ar @{=}[r]&
 W  
}
\]
where $\lam\in\mu_{A}(\WW)$ and $\orb$ is the coadjoint orbit passing
through $\lam$. Therefore it is a sensible generalization of the smooth
case to interpret the reduced space
$\mu_{A}^{-1}(\orb)/G = \WW\spr{\mathcal{O}}G$ as a typical
stratified symplectic leaf of the stratified Poisson space $W$. 
The following thus generalizes the result of Marsden and Perlmutter
\cite[Theorem 4.3]{MP00} to the case of a non-free but single orbit
type action of $G$ on $Q$.

Let $\orb$ be a
coadjoint orbit in the image of the \momap $\mu_{A}: \WW\to\gu^{*}$,
and let $(L)$ be in the isotropy lattice of the $G$-action on
$\WW$ such that
$\lorb{\WW}:=\mu_{A}^{-1}(\orb)\cap\WW_{(L)}\neq\emptyset$. 
Then we define
\[
 \lorb{\iota}: \lorb{\WW}\hookto\WW,
\]
the canonical embedding, and the orbit projection mapping
\[
 \lorb{\pi}: \lorb{\WW}\toto\lorb{\WW}/G=:(\WW\sporb G)_{(L)}.
\]
Further, we denote the Kirillov-Kostant-Souriau symplectic form on
$\orb$ by $\omkks$, that is
$\omkks(\lam)(\ad^{*}(X).\lam,\ad^{*}(Y).\lam) = \vv<\lam,[X,Y]>$. 
Remember from Proposition \ref{prop:WW} that the symplectic structure
on $\WW=\WW(A)$ is denoted by $\sigma = \psi^*\Om$.

\begin{theorem}[Gauged symplectic reduction]\label{thm:WWred}
Let $Q=Q_{(H)}$
and let $A$ be a generalized connection form on $\pi: Q\toto Q/G$.
Let $\orb$ be a
coadjoint orbit in the image of the \momap $\mu_{A}: \WW\to\gu^{*}$,
and let $(L)$ be in the isotropy lattice of the $G$-action on
$\WW$ such that
$\lorb{\WW}:=\mu_{A}^{-1}(\orb)\cap\WW_{(L)}\neq\emptyset$. Then the
following are true.
\begin{enumerate}[\up (1)]
\item 
The smooth manifold $(\WW\sporb G)_{(L)}$ is a typical symplectic
stratum of the singular symplectic space $\WW\sporb G$. The smooth
symplectic manifold 
  \[
   (\orb\spr{0}H)_{(L_0)_H^G} := (\orb\cap\ann\ho)_{(L_0)_H^G}/H,
  \]
where $L_0$ is an isotropy subgroup of
the $H$-action on $\orb\cap\ann\ho$,
is a disjoint union of 
typical smooth symplectic strata of the singular symplectic space
$\orb\spr{0}H$. 
\item 
The symplectic stratum $(\WW\sporb G)_{(L)}$ can be globally
described as
\[ 
  (\WW\sporb G)_{(L)}
  =
  T^{*}(Q/G)\times_{Q/G}(\bscorb)_{(L)}/G
\]
whence it is the total space of the smooth symplectic fiber bundle
\[
\xymatrix{
 {(\orb\spr{0}H)_{(L_0)_H^G}}
 \ar @{^{(}->}[r]&
 {(\WW\sporb G)_{(L)}}
 \ar[r]^{\chi}
 &{T^{*}(Q/G).}
}
\]
Hereby $L_{0}\subset H$ is an
isotropy subgroup of the induced $H$-action on $\orb$ which is
conjugate to $L$ within $G$. 
\item 
The symplectic structure $\lorb{\sigma}$ on $(\WW\sporb G)_{(L)}$
is uniquely determined and given by the formula
\[
 (\lorb{\pi})^{*}\lorb{\sigma}
 =
 (\lorb{\iota})^{*}\sigma - (\mu_{A}|\lorb{\WW})^{*}\omkks.
\]
Therefore, 
\[
 \lorb{\sigma} = \chi^*\Om^{Q/G} - \lorb{\beta}
\]
where $\lorb{\beta}\in\Om^{2}((\WW\sporb G)_{(L)})$ is defined by 
\[
 (\lorb{\pi})^{*}\lorb{\beta}
 =
  (\lorb{\iota})^{*}dB + (\mu_A|\lorb{\WW})^{*}\omkks.
\]
Finally, $B$ is the form that was introduced in Proposition
\ref{prop:WW}.
\item  
The stratified symplectic space can be globally described as 
\[ 
 \WW\sporb G
 = 
 T^{*}(Q/G)\times_{Q/G}\bscorb/G
\]
whence it is the total space of
\[
\xymatrix{
 {\orb\spr{0}H}
 \ar @{^{(}->}[r]&
 {\WW\sporb G}
 \ar[r]
 &{T^{*}(Q/G)}
}
\]
which
is a stratified symplectic fiber bundle with singularities
confined to the fiber direction.
\end{enumerate}
\end{theorem}

\begin{proof}
Assertion (1). 
This is an implication of the general theory of stratified
symplectic reduction -- see Ortega and Ratiu~\cite[Section 8.4]{OR04}
or Section~\ref{s:sg_com_red} for a statement of the relevant theorem
and Section~\ref{sec:prel} for the notation.

To see that 
$
 (\orb\spr{0}H)_{(L_0)_H^G}
$
is a union of 
typical smooth symplectic strata of the singular symplectic space
$\orb\spr{0}H$ note firstly that 
$
 (\orb\spr{0}H)_{(L_0)_H^G}
$
is,
according to the proof of Lemma~\ref{lem:weirdo-type},
a typical stratum of $(T^*G)\spr{0}H\sporb G$. 
By the corollary of
Lemma~\ref{lem:actions} (with $S=\set{\textup{point}}$ and $\beta = 0$) 
there is an isomorphism 
$
 (T^*G)\spr{0}H\sporb G 
 \cong (T^*G)\sporb G\spr{0}H 
 = 
 \orb\spr{0}H
$, 
$
 [[(g,\lam)]_H]_G
 \mapsto  
 [[(g,\lam)]_G]_H
 =
 [\lam]_H$
of singular symplectic spaces, whence strata are mapped
symplectomorphically onto unions of strata.

Assertion (2). 
The description of the stratum $(\WW\sporb G)_{(L)}$ follows from
Proposition~\ref{prop:WW}.

We know from (1) that all spaces involved in the
diagram really are smooth. As in the proof of Lemma \ref{lm:WWorb}
let $q_{0}\in Q$ with $G_{q_{0}}=H$, $S$ a slice at $q_{0}$, and
$U\cong G/H\times S$ a tube around the orbit $G.q_{0}$. Then we get
the local description
\begin{align*}
 (\WW\sporb G)_{(L)}|U
 &=
 T^{*}S\times_{S}(\bsc_{q\in U}\orb\cap\ann\gu_{q})_{(L)}/G\\
 &\cong
 T^{*}S\times_{S}S\times(\orb\cap\ann\ho)_{(L_0)_H^G}/H\\
 &=
 T^{*}S\times(\orb\spr{0}H)_{(L_0)_H^G}
\end{align*}
as claimed. 
The bundle is symplectic by Theorem~\ref{thm:bun_pic}. 

Assertion (3). The defining property of the reduced symplectic form
$\lorb{\sigma}$, namely,
\[
 (\lorb{\pi})^{*}\lorb{\sigma}
 =
 (\lorb{\iota})^{*}\sigma - (\mu_{A}|\lorb{\WW})^{*}\omkks
\]
is a well-established fact, 
see e.g.\ Bates and Lerman~\cite[Proposition~11]{BL97}. 
Thus it is clear from Proposition~\ref{prop:WW} that 
\[
 \lorb{\sigma}
 =
 \chi^*\Om^{Q/G}-\lorb{\beta}
\]
and it remains to check that 
$\lorb{\beta}$ is a well defined two-form on $(\WW\sporb G)_{(L)}$.
To see this notice firstly that  
\[
 \tilde{\beta} 
 :=
 (\lorb{\pi})^{*}\lorb{\beta}
 =
  (\lorb{\iota})^{*}dB + (\mu_A|\lorb{\WW})^{*}\omkks
 \in \Om^2(\lorb{\WW})
\]
is
$G$-invariant. Furthermore, we claim that $\tilde{\beta}$ is
horizontal, i.e., vanishes upon insertion of a vertical
vector field. Indeed, let 
$
 (q,\eta,\lam)\in\lorb{\WW} 
$
and consider $\zeta_{Z_1}(q,\eta,\lam),\xi_2\in
T_{(q,\eta,\lam)}\lorb{\WW}$.
We proceed locally as in the proof of Proposition~\ref{prop:WW} so
that there is a splitting of tangent vectors as 
$\zeta_{Z_1}(q,\eta,\lam) = (\zeta_{Z_1}(q),0,\ad^*(Z_1).\lam)$ and
$\xi_2 = (q_2',\eta_2',\ad^*(Y).\lam)$. 
Therefore,
\begin{align*}
 \tilde{\beta}_{(q,\lam)}(&\zeta_{Z_1}(q,\eta,\lam),\xi_2)\\
  &=
  \vv<\ad^{*}(Z_{1}).\lam,Z_{2}>
  - \vv<\ad^{*}(Y).\lam,Z_{1}>
  + 0
  + \vv<\lam,[Z_{1},Z_{2}]>
  + \vv<\lam,[Z_{1},Y]>\\
  &= 
    - \vv<\lam,[Z_{1},Z_{2}]>
    - \vv<\lam,[Z_{1},Y]>
    + \vv<\lam,[Z_{1},Z_{2}]>
    + \vv<\lam,[Z_{1},Y]> 
  = 0.
\end{align*}   
That is, $\tilde{\beta}$ is a basic form and thus descends to a form
$\lorb{\beta}$. 

Assertion (4) is a pasting together of the results in (2). 
\end{proof}

\begin{corollary}\label{cor:WWred}
Let $\orb$ be a
coadjoint orbit in the image of the \momap $\mu_{A}: \WW\to\gu^{*}$,
and let $(L)$ be in the isotropy lattice of the $G$-action on
$\WW$ such that
$\lorb{\WW}:=\mu_{A}^{-1}(\orb)\cap\WW_{(L)}\neq\emptyset$. Assume
further that there is a global slice $S$ such that $Q\cong G/H\times
S$. Then we have the global description 
\[
 \WW\sporb G =
 T^{*}S\times\orb\spr{0}H. 
\]
Moreover, the reduced symplectic form
$\lorb{\sigma}$ on a symplectic stratum 
$(\WW\sporb G)_{(L)} = T^{*}S\times(\orb\spr{0}H)_{(L_0)_H^G}$ 
is given by the formula 
\[
 \lorb{\sigma} = \Om^{Q/G}-\Om^{\mathcal{O}}_{(L_0)_H^G}
\]
where $\Om^{\mathcal{O}}_{(L_0)_H^G}$ is the canonically reduced 
symplectic form on
$(\orb\spr{0}H)_{(L_0)_H^G}$, and $L_{0}$ is a subgroup of $H$ which
is conjugate to $L$ within $G$.
\end{corollary}

\begin{proof}
This is an immediate consequence of Theorems \ref{thm:bun_pic} and
\ref{thm:WWred}.
\end{proof}

\section{Spin Calogero-Moser systems}\label{sec:cms}

In this section we give a mechanical application of
Theorem~\ref{thm:WWred} to obtain spin Calogero-Moser models. 
This approach follows, in essence, the idea of Kazhdan, Kostant and
Sternberg~\cite{KKS78} that such models may be obtained via
projection of geodesic systems on Lie groups or Lie algebras.
This construction can be
carried out in various guises and at different levels of generality.
For the simplest case 
we describe the way Theorem~\ref{thm:WWred} allows to understand this
projection procedure in the following subsection. The emphasis is here
on the use of the mechanical connection.

\subsection{The construction based on cotangent bundle reduction}\label{sub:cms-constr}
Let $G$ be a (real or complex)
simple Lie group, $\gu$ its Lie algebra, $\ho$ a Cartan sub-algebra,
and $H$ a corresponding Cartan subgroup. 
Then we consider either $Q = G_{(H)}$ acted upon by $G$ via
conjugation or $Q = \gu_{(H)}$ acted upon by $\Ad(G)$.
Note that $G_{(H)}$ is open dense in $G$ and $\gu_{(H)}$ is open dense
in $\gu$.  
We will see below that choosing $Q=\gu_{(H)}$ leads to Calogero-Moser
systems with rational potential while choosing $Q=G_{(H)}$ leads to
Calogero-Moser systems with trigonometric potential.
As in the
construction of the previous sections we may then consider the lifted
$G$-action on $T^*Q$ which is Hamiltonian with equivariant \momap
$\mu: T^*Q\to\gu^*$. Since $G$ is assumed simple we can use the
Killing form $B$
to identify $T^*Q\cong_B TQ$ and $\gu^*\cong_B\gu$ whence $\mu$ becomes a
mapping $\mu: TQ\to\gu$. Let $\orb$ be an adjoint orbit in the image
of $\mu$. 
Via the Killing form we may thus define a mechanical connection $A$ on
$Q\toto Q/G$ as in Subsection~\ref{sub:mech-con}.
By Theorem~\ref{thm:WWred} and its corollary the singular symplectic
quotient of $TQ$ at $\orb$ is therefore given by
\[
 \mu^{-1}(\orb)/G
 =
 TQ\sporb G
 \cong_{A}
 T(Q/G)\times\orb\spr{0}H,
\]
and this space is called a \caps{spin Calogero-Moser space}. 
(Recall that $\orb\spr{0}H = (\orb\cap\ho^{\bot})/H$.)
This terminology is justified as follows. Use the left multiplication
in the group to trivialize the tangent bundle as $TQ =
Q\times\gu$. Let $\mathcal{H}: Q\times\gu\to\R$,
$(q,X)\mapsto\by{1}{2}B(X,X)$ denote the free Hamiltonian. In the
notation of Proposition~\ref{prop:WW} the reduced Hamiltonian system 
on $T(Q/G)\times\orb\spr{0}H$ is thus given by Hamiltonian reduction
of $(\WW,\sigma,\psi^*\mathcal{H})$ at $\orb$. 
Since the Hamiltonian in this
picture is given by 
\[
 (\psi^*\mathcal{H})(q,\eta,\lam)
 =
 \by{1}{2}B(\eta,\eta) + \by{1}{2}B(A_q^*(\lam),A_q^*(\lam))
\]
for $(q,\eta,\lam)\in\WW$ one thus needs to compute $A_q^*(\lam)$. In
fact, it obviously suffices to compute $A_q^*(\lam)$ for $q\in Q_H$. 
The crucial point now is that with respect to the identifications
$T^*Q = TQ$ and $\gu^* = \gu$ we have
\[
 A_q^*(\lam) = \zeta(q)(\check{\ine}_q^{-1}(\lam))
\]
where $\check{\ine}_q: \gu_q^{\bot}\to\gu_q^{\bot}$
is the $G_q$-equivariant \iso obtained from the locked inertia
tensor $\ine$. 
This object can be computed using structure theory, and we shall do so
in the next paragraph.

\subsubsection{The rational case}
Assume that $\gu$ is a complex simple Lie algebra. 
(The case of real simple Lie algebras works analogously.)
Let $Q=\gu_{(H)}$, let $\Delta\subset\ho^*$ be a root system,
$\Delta^+$ a system of positive roots, and $\gu
= \ho\oplus\oplus_{\alpha\in\Delta}\gu_{\alpha}$ the corresponding
root space decomposition. For each $\alpha\in\Delta$ we choose a
vector $E_{\alpha}\in\gu_{\alpha}$ such that $B(E_{\alpha},E_{-\alpha})
= 1$. 
Assume that $q\in Q_H = \ho_{\textup{reg}}$ where $\ho_{\textup{reg}}$
denotes the set of regular elements.
Then we may write 
$\lam\in\orb\cap\ho^{\bot}$ as 
$\lam =
\sum_{\alpha\in\Delta}\lam_{\alpha}E_{\alpha}$ 
and 
$
 \smash{
 Z :=
 \check{\ine}_q^{-1}(\lam)\in\ho^{\bot}
 }
$ 
as 
$Z = \sum_{\alpha\in\Delta}z_{\alpha}E_{\alpha}$. 
Since $\zeta_Z(q) = \ad(Z).q$
it follows that
\begin{align*}
 \lam_{\alpha}
 &=  
 B(\lam,E_{-\alpha})
 =
 \ine_q(Z,E_{-\alpha})
 =
 B(\zeta_Z(q),\zeta_{E_{-\alpha}}(q))
 =
 B(Z,-\alpha(q)^2 E_{-\alpha})\\
 &=
 -z_{\alpha}\alpha(q)^2.
\end{align*}
Therefore, we find that
$
 \smash{
 Z 
 = 
 -\sum_{\alpha\in\Delta}\lam_{\alpha}\alpha(q)^{-2}E_{\alpha}
 }
$,
and 
\[
\tag{R}\label{equ:R}
 A_q^*(\lam)
 =
 \zeta_Z(q)
 =
 +\ad(q)\sum_{\alpha\in\Delta}\by{\lam_{\alpha}}{\alpha(q)^2}E_{\alpha}
 =
 \sum_{\alpha\in\Delta}\by{\lam_{\alpha}}{\alpha(q)}E_{\alpha}.
\]
This in turn implies that the reduced Hamiltonian
$(\psi^*\mathcal{H})_0$ is given by
\begin{align*}
 (\psi^*\mathcal{H})_0(q,\eta,[\lam]_H)
 &=
 \by{1}{2}B(\eta,\eta) + 
  \by{1}{2}B(\sum_{\alpha\in\Delta}\by{\lam_{\alpha}}{\alpha(q)}E_{\alpha}
             ,\sum_{\alpha\in\Delta}\by{\lam_{\alpha}}{\alpha(q)}E_{\alpha})\\
 &=
 \by{1}{2}B(\eta,\eta) 
 - 
 \sum_{\alpha\in\Delta^+}\by{\lam_{\alpha}\lam_{-\alpha}}{\alpha(q)^2}
\end{align*}
where $(q,\eta,[\lam]_H)\in T(Q/G)\times\orb\spr{0}H =
C_{\textup{reg}}\times\ho\times\orb\spr{0}H$ and where $C_{\textup{reg}}$ denotes
the interior of a Weyl chamber. 
This function is the Hamiltonian of the rational
Calogero-Moser system with spin.

\subsubsection{The trigonometric case}
Let $Q=G_{(H)}$ and continue the notation regarding the structure
theoretic objects of the previous paragraph. Assume $q = \exp{a}\in Q_H =
H_{\textup{reg}}$ whence $a\in\ho_{\textup{reg}}$. Consider
$\lam\in\orb\cap\ho^{\bot}$ and $Z$ as above. Since $\zeta_Z(q) =
\Ad(q^{-1}).Z-Z = (e^{-\textup{ad}(a)}-1).Z$ we may compute
\[
 \lam_{\alpha}
 =
 B(\lam,E_{-\alpha})
 =
 \ine_q(Z,E_{-\alpha})
 =
 (e^{-\alpha(a)}-1)(e^{\alpha(a)}-1)z_{\alpha}
 =
 (2-2\cosh \alpha(a))z_{\alpha}
\]
and thus obtain the adjoint to the mechanical connection
\[
\tag{T}\label{equ:T}
 A_q^*(\lam)
 =
  \zeta_Z(q)
 =
   \by{1}{2}\sum_{\alpha\in\Delta}\by{2}{e^{-\alpha(a)}-1}\lam_{\alpha}E_{\alpha}
 =
  \by{1}{2}\sum_{\alpha\in\Delta}\lam_{\alpha}E_{\alpha}
  +
  \by{1}{2}\sum_{\alpha\in\Delta}\coth(\by{-\alpha(a)}{2})\lam_{\alpha}E_{\alpha}
\]
in the same way as above.
Therefore, the reduced Hamiltonian
$(\psi^*\mathcal{H})_0$ is given by
\[
 (\psi^*\mathcal{H})_0(q,\eta,[\lam]_H)
 =
  \by{1}{2}B(\eta,\eta) + \by{1}{2}B(\zeta_Z(q),\zeta_Z(q))
 =
  \by{1}{2}B(\eta,\eta) 
   -
  \by{1}{4}\sum_{\alpha\in\Delta^+}\by{\lam_{\alpha}\lam_{-\alpha}}{\sinh^2\alpha(a)}
\]
where $(q,\eta)\in T(Q_{(H)}/G) = T(H_{\textup{reg}}/W)$ and
$[\lam]_H\in\orb\spr{0}H$, and where $W=N(H)/H$ is the Weyl group.  
This function is the Hamiltonian of the trigonometric
Calogero-Moser system with spin.

\begin{rem}[Mechanical connection \& classical dynamical
  $r$-matrix]\label{rem:r+A}
It is noted in the introduction that the idea of obtaining
Calogero-Moser systems through Hamiltonian reduction is originally due
to Kazhdan, Kostant and Sternberg~\cite{KKS78}.
However, a completely different approach to obtain such systems was
taken by Li and Xu~\cite{LX00,LX02}. They used 
an analysis based on
classical dynamical
$r$-matrices (associated to complex simple Lie algebras) 
to directly write down 
the Hamiltonian of spin Calogero-Moser systems (associated to complex simple
Lie algebras). 
(See also Feh\'{e}r and Pusztai~\cite[Section~2]{FP06} for an outline of
this construction.)
This approach is based on the classification of classical
dynamical $r$-matrices of Etingof and Varchenko~\cite{EV98}.
It was noticed by Feh\'{e}r and Pusztai~\cite{FP06} that 
one can obtain 
the same
Calogero-Moser models which appear in \cite{LX00,LX02} through
Hamiltonian reduction of cotangent bundles. 
Moreover, 
constructing certain new spin Calogero-Moser models
it was observed by \cite[Proof of Prop.~3]{FP06} that 
dynamical $r$-matrices appear in the process of Hamiltonian reduction
of 
$T^*G$ where $G$ is a complex or real simple Lie group acting on itself
by twisted conjugation. 
However, the relationship
of the $r$-matrix and the reduction approach was still mysterious
in the sense that there was no explanation for it other than the
computations obviously yielding the correct results.
We claim that this
relationship can be further 
explained in a geometric framework using the
mechanical connection.\footnote{This point of view essentially evolved
  during discussions with Laszlo Feh\'{e}r.} 
Indeed,
in \cite{EV98} a classical dynamical $r$-matrix 
associated to a complex simple Lie algebra $\gu$ 
is defined as a
meromorphic function $r: \ho^*=_B\ho\to\gu\otimes\gu$ which satisfies 
the classical dynamical Yang-Baxter equation 
(CDYBe) and certain other conditions (\cite[Subsection~3.2]{EV98})
in the completed tensor product
$\gu\hat{\otimes}\gu$. Via the \iso $\gu\otimes\gu \cong
\hom(\gu^*,\gu) =_B \hom(\gu,\gu)$ we may think of a classical
dynamical $r$-matrix as a meromorphic function 
$R:
\ho\to\hom(\gu,\gu)$ subject to the appropriate equations. 

For the \emph{rational case} let us use Equation~(\ref{equ:R}) to define
a holomorphic function $R: \ho_{\textup{reg}}\to\hom(\gu,\gu)$ by
\[
 R(q)(\lam)
 :=
 A_q^*(\lam)
 = 
 \sum_{\alpha\in\Delta}\by{\lam_{\alpha}}{\alpha(q)}E_{\alpha}
\]
By \cite[Theorem~3.2]{EV98} any classical dynamical $r$-matrix 
associated to a complex simple Lie algebra $\gu$ with
coupling constant $\epsilon=0$ is of this form 
(where $X=\Delta$, $C=0$ and
$\nu=0$ in the
notation of \cite[Theorem~3.2]{EV98}).

For the \emph{trigonometric case} let us use Equation~(\ref{equ:T}) 
 to define
a holomorphic function $R: \ho_{\textup{reg}}\to\hom(\gu,\gu)$ by
$R(a)(\lam)=\by{1}{2}\lam$ for $\lam\in\ho$, and 
\begin{align*}
 R(a)(\lam)
 &:=
  A_{\exp a}^*(\lam)
 = 
  \by{1}{2}\sum_{\alpha\in\Delta}\lam_{\alpha}E_{\alpha}
   +
  \by{1}{2}\sum_{\alpha\in\Delta}\coth(\by{-\alpha(a)}{2})\lam_{\alpha}E_{\alpha}\\
 &=
  i_{\lam}
  (\by{1}{2}\Om
   +
  \by{1}{2}\sum_{\alpha\in\Delta}\coth(\by{\alpha(a)}{2})E_{\alpha}\otimes E_{-\alpha})
\end{align*}
for $\lam\in\ho^{\bot}$.
Here 
$\Om 
 = \sum_{i=1}^{l}x_i\check{B}(x_i) 
   + 
   \sum_{\alpha\in\Delta}E_{\alpha}\check{B}(E_{-\alpha})$ 
is the Casimir element of $(\gu,B)$
where $x_1,\ldots,x_l$ is an orthonormal basis of $\ho$,
and 
$X=\Delta_+$,
$C_{i,j}=0$,
$\epsilon = 1$ and $\nu=0$
in the notation of \cite[Theorem~3.10]{EV98}.
By \cite[Theorem~3.10]{EV98} any classical dynamical $r$-matrix 
associated to a simple Lie algebra with
coupling constant $\epsilon=1$ is of this form. 

We view this as a new geometric explanation of why it is
possible to associate Calogero-Moser systems to classical dynamical
$r$-matrices. It remains a goal for future work 
to find a general
relationship between the condition (i.e., CDYBe) defining classical
dynamical $r$-matrices and the properties 
(such as $\curv^A = 0$)
of the mechanical connection.  
\end{rem}

\subsection{$\SL(m,\C)$ by hand}\label{sub:kks}
As an example consider $G=\SL(m,\C)$. Here we work along the lines of \kks
\cite[Section 2]{KKS78} who considered the case $G=\SU(m,\C)$. See
also \aklm \cite[Section 5.7]{AKLM03}. 
The point to this example is that we try to say as much as possible
about the reduced phase space by using an \emph{ad hoc} approach.

Let $\orb=\Ad(G)Z_{0}$ be an
orbit passing through a semi-simple element $Z_{0}$. Consider $(a,\alpha)\in
G_{r}\times\gu$ with $\alpha-a\alpha a^{-1} = \mu(a,\alpha) = Z$. 
Note that $\mu: G_r\times\gu\to\gu$ is the equivariant \momap of the
action which is obtained by lifting the $G$-action on $G_r$ by
conjugation 
to the (co-)tangent bundle $G_r\times\gu$.  
As usual $G_{r}$ denotes the set of regular elements, that is, $G_{r}$
consists of those matrices that have $m$ different
eigenvalues. Moreover, we let $H$ denote the subgroup of diagonal
matrices, and $H_{r} := H\cap G_{r}$. 
Via the $\Ad(G)$-action we can bring $a$ in diagonal form with entries
$a_{i}\neq a_{j}$ for $i\neq j$.
Since $Z_{ij} =
\alpha_{ij}-\by{a_{i}}{a_{j}}\alpha_{ij}$ the following are
coordinates on 
$(\mu^{-1}(\orb)\cap(G_{r}\times\gu))/\Ad(G)
 = T^*G_r\sporb G$:
\begin{itemize}
\item $a_{i}$ for $i=1,\ldots,m$.
\item $\alpha_{i}:=\alpha_{ii}$ for $i=1,\ldots,m$. 
\item $\alpha_{ij}=(1-\by{a_{i}}{a_{j}})^{-1}Z_{ij}$ for $i\neq j$.
\end{itemize}
These coordinates give an identification
\[
 (\mu^{-1}(\orb)\cap(G_{r}\times\gu))/\Ad(G) =
   (T^{*}H_{r}\times(\orb\cap\ho^{\bot})/\Ad(H))/W
\]
where $W=N(H)/H$ is the Weyl group.
\caps{Claim:} \emph{If $\orb$ is an orbit which is of minimal non-zero
  dimension then we have that $\orb\cap\ho^{\bot}/\Ad(H) =
  \set{\textup{point}}$. Moreover, the reduced phase space can be
  described as
  $(\mu^{-1}(\orb)\cap(G_{d}\times\gu))/\Ad(G) \cong T^{*}H_{r}/W$.}

Here $G_{d}$ denotes the open and dense subset of all diagonable elements in
$\SL(m,\C)$. Indeed, let $\mu(a,\alpha) = Z\in\orb\cap\ho^{\bot}$
with $a$ in diagonal form. Thus
$Z=vw^{t}-cI$ where $c:=\by{1}{m}\vv<v,w>\neq0$, $v,w\in\C^{m}$,
and $w^{t}$ is the transposed to the column vector $w$. Since
$Z\in\ho^{\bot}$ we infer that $v_{i}w_{i}=c$. Hence 
\[
 \orb\cap\ho^{\bot} =
  \set{(\by{c}{v_{1}}v,\dots,\by{c}{v_{m}}v)-cI:
         v_{i}\in\C\setminus\set{0}}. 
\]
Take such an
$(\by{c}{v_{1}}v,\dots,\by{c}{v_{m}}v)-cI=:Z_{1}$. Let
$h=\prod_{i=1}^{m}v_{i}\cdot\diag{v_{1}^{-1},\dots,v_{m}^{-1}}$. Then we
can bring $Z_{1}$ into the normal form
$\Ad(h)Z_{1} = c(1)_{ij}-cI$ where $(1)_{ij}$ denotes the $m\times
m$-matrix with all entries equal to $1$. Finally note that
$\alpha_{ij}-\by{a_{i}}{a_{j}}\alpha_{ij} = \by{c}{v_{j}}v_{i} \neq 0$
implies that $a=\diag{a_{1},\dots,a_{m}}$ is actually regular. 

The coordinates for $T^*G_r\sporb G$ found above by evaluating the momentum constraint
equation and factoring out the $G$-action have been the motivating
point for the formulation of the general Theorem~\ref{thm:WWred}.

\subsection{Application: Hermitian matrices}
Consider $V$ the space of complex Hermitian $n\times n$ matrices as
the configuration space to start from. This space shall be endowed
with the inner product $V\times V\to\R$,
$(a,b)\mapsto\tr(ab)$. Moreover, we let $G=\SU(n,\C)$ act on $V$ by
conjugation. Clearly this action leaves the trace form invariant. Via
the inner product we can trivialize the cotangent bundle as $T^{*}V =
V\times V^{*} = V\times V$, and the cotangent lifted action of $G$ is
simply given by the diagonal action. 
The canonical symplectic form on $T^{*}V$ is given by
\[
 \Om_{(a,\alpha)}((a_{1},\alpha_{1}),(a_{2},\alpha_{2})) =
 \tr(\alpha_{2}a_{1})-\tr(\alpha_{1}a_{2}). 
\]
The free Hamiltonian on $T^{*}V=V\times V$ is given by 
\[
 \hfree: (a,\alpha)\longmapsto\by{1}{2}\tr(\alpha\alpha).
\]
Trajectories of this Hamiltonian are given by straight lines of the
form $t\mapsto a+t\alpha$ in the configuration space $V$. 

Let us further identify
$\su(n)^{*}=\su(n)$ via the Killing form. The momentum mapping is then
given by 
\[
 \mu: (a,\alpha)\longmapsto[a,\alpha] = \ad(a).\alpha.
\]
Consider also an orbit
$\orb$ together with its canonically induced symplectic structure in
the image of the momentum mapping. 

\textbf{Assumption:} The orbit $\orb$ is such that
$\mu^{-1}(\orb)\subeq V_{r}\times V$. Here $V_{r}$ denotes the set of
regular elements in $V$ with respect to the $G$ action. 

This assumption is, for example, fulfilled if the projection from $\orb$
to any root space is non-trivial.
On the other hand, if the assumption is not satisfied for a particular
orbit $\orb$ one can also consider the restricted $G$-action on $V_r$
and proceed with reduction of the Hamiltonian system
$(T^*V_r,\Om,\hfree)$ at the orbit level $\orb$. This has, however,
the disadvantage that the Hamiltonian flow lines may leave
$T^*V_r\sporb G$ in finite time. 

Let $\Sigma$ denote the subspace of $V$ consisting of diagonal
matrices. Then $\Sigma$ is a section of the $G$-action on $V$, see
Section \ref{sec:pol_rep}. Further, we define
$\Sigma_{r}:=V_{r}\cap\Sigma$. Within $\Sigma$ we choose the positive Weyl
chamber $C:=\set{\diag{q_{1},\dots,q_{n}}: q_{1}>\ldots>q_{n}}$ so
that $C = \Sigma/W = V/G$ 
where 
$W = W(\Sigma) = N_{G}(\Sigma)/Z_{G}(\Sigma)$. 
Thus $C_{r}:=\Sigma_{r}\cap C$ 
may be considered as
a global slice for the $G$-action on $V_{r}$ so that 
$G/M\times C_{r}\cong V_{r}$, $(gM,a)\mapsto g.a$ 
where $M:=Z_{G}(\Sigma_{r})=Z_{G}(\Sigma)$. That
is, $M$ is the subgroup of $\SU(n)$ consisting of diagonal matrices
only. Now we may apply Corollary \ref{cor:WWred} to get 
\[
 T^{*}V\sporb G
 =
 T^{*}C_{r}\times\orb\spr{0}M
\]
as symplectic stratified spaces. The strata are of the form 
\[
 (T^{*}V\sporb G)_{(L)}
 =
 T^{*}C_{r}\times(\orb\spr{0}M)_{(L_0)_M^G}
\]
where $L_{0}$ is a subgroup of $M$ conjugate to $L$ within $G$. 
Moreover, the reduced symplectic structure $\lorb{\sigma}$ on
$(T^{*}V\sporb G)_{(L)}$ is of product form, i.e., 
\[ 
 \lorb{\sigma} = 
  \Om^{C_{r}}-\Om^{\mathcal{O}}_{(L_0)_M^G}
\]
where $\Om^{\mathcal{O}}_{(L_0)_M^G}$ is the canonically reduced
symplectic form on $(\orb\spr{0}M)_{(L_0)_M^G}$.  

From the general theory (Theorem~\ref{thm:sing_spr})
we know that the Hamiltonian $\hfree$ reduces
to a Hamiltonian $\hcml$ on the stratum $(T^ {*}V\sporb G)_{(L)}$, and
that integral curves of $\hfree$ project to integral curves of
$\hcml$. 
In particular, the dynamics remain confined to the symplectic stratum. 
The reduced Hamiltonian is thus given by 
\[
 \hcml(q,p,[\lam]) = \hfree(q,p+A_{q}^{*}(\lam))
\]
where $[\lam]$ is the class of $\lam$ in $(\orb\spr{0}M)_{(L_0)_M^G}$ and
$A_{q}^{*}: \gu_{q}^{\bot}=\mo^{\bot}\to T_{q}(G.q)=\Sigma^{\bot}$ 
is the point wise dual to the mechanical connection as
introduced in Section \ref{sec:gauged_red}. Assume that
$q=\diag{q_{1}>\ldots>q_{n}}$ and that
$\lam=(\lam_{ij})_{ij}\in(\orb\cap\mo^{\bot})_{(L_0)_M^G}$. Then 
\[
 A_{q}^{*}(\lam)_{ij} = \by{\lam_{ij}}{q_{i}-q_{j}}
 \textup{  for  } i\neq j,
 \textup{  and  }
 A_{q}^{*}(\lam)_{ii} = 0.
\]
Therefore, for $p=\diag{p_{1},\ldots,p_{n}}\in\Sigma$ and $q,[\lam]$
as introduced we obtain
\begin{align*}
 \hcml(q,p,[\lam])
 &= 
 \by{1}{2}\tr(p)^{2}+\by{1}{2}\tr(A_{q}^{*}(\lam))^{2}
 =
 \by{1}{2}\sum_{i=1}^{n}p_{i}^{2}
 + \by{1}{2}\sum_{i\neq
      j}\by{\lam_{ij}\lam_{ji}}{(q_{i}-q_{j})(q_{j}-q_{i})}\\
 &=
 \by{1}{2}\sum_{i=1}^{n}p_{i}^{2}
 + \sum_{i>j}\by{|\lam_{ij}|^{2}}{(q_{i}-q_{j})^{2}}
\end{align*}
since $\lam_{ji}=-\overline{\lam_{ij}}$ and
$\tr(pA_{q}^{*}(\lam))=\tr(A_{q}^{*}(\lam)p)=0$. 
This is the Hamiltonian
function of the Calogero-Moser system with spin. Integrability of this
system in the non-commutative sense is proved in the next section in a
more general context.

\subsection{Application: Polar representations of compact Lie groups}\label{sub:pol-cm}

The idea of considering polar representations of compact Lie groups to
obtain new versions of Spin Calogero-Moser systems is due to \aklm
\cite{AKLM03}.
 
As in Section  \ref{sec:pol_rep} let $V$ be a real Euclidean
vector space and $G$ a connected compact Lie group that acts on $V$ via
a polar representation. Via the inner product we consider the
cotangent bundle of $V$ as a product $T^{*}V=V\times V$. The canonical
symplectic form $\Om$ is thus given by 
\[
 \Om_{(a,\alpha)}((a_{1},\alpha_{1}),(a_{2},\alpha_{2})) =
 \vv<\alpha_{2},a_{1}>-\vv<\alpha_{1},a_{2}>
\]
where $\vv<\phantom{a},\phantom{a}>$ is the inner product on $V$. 
The free Hamiltonian on $T^{*}V=V\times V$ is given by 
\[
 \hfree: (a,\alpha)\longmapsto\by{1}{2}\vv<\alpha,\alpha>.
\]
Trajectories of this Hamiltonian are given by straight lines of the
form $t\mapsto a+t\alpha$ in the configuration space $V$. 

Of course, the
cotangent lifted action of $G$ is just the diagonal action
of $G$ on $V\times V$. By Section \ref{sec:pol_rep} we may 
think of the action by $G$ on $V$ as a symmetric space representation
and thus consider
$\gu\oplus V =: \lo$ as a real semi-simple Lie algebra with Cartan decomposition
into $\gu$ and $V$, and with bracket relations
$[\gu,\gu]\subeq\gu$, $[\gu,V]\subeq V$, and $[V,V]\subeq V$. The
momentum mapping corresponding to the $G$-action on $T^{*}V=V\times V$
is now given by 
\[
 \mu: V\times V\longto\gu^{*}=\gu,
 \quad
 (a,\alpha)\longmapsto[a,\alpha]=\ad(a).\alpha
\]
where we identify $\gu=\gu^{*}$ via an $\Ad(G)$-invariant inner product. 
Consider also an orbit
$\orb$ together with its canonically induced symplectic structure in
the image of the momentum mapping. 

\textbf{Assumption:} The orbit $\orb$ is such that
$\mu^{-1}(\orb)\subeq V_{r}\times V$. Here $V_{r}$ denotes the set of
regular elements in $V$ with respect to the $G$ action. 
 
We proceed as above, and let $\Sigma$ denote a fixed section of the
$G$-action on $V$, consider $C$ a Weyl chamber in $\Sigma$, and put
$M:=Z_{G}(\Sigma)$. 
We may apply Corollary \ref{cor:WWred} to get 
\[
 T^{*}V\sporb G
 =
 T^{*}C_{r}\times\orb\spr{0}M
\]
as symplectic stratified spaces. The strata are of the form 
\[
 (T^{*}V\sporb G)_{(L)}
 =
 T^{*}C_{r}\times(\orb\spr{0}M)_{(L_0)_M^G}
\]
where $L_{0}$ is a subgroup of $M$ conjugate to $L$ within $G$. 
Moreover, the reduced symplectic structure $\lorb{\sigma}$ on
$(T^{*}V\sporb G)_{(L)}$ is of product form, i.e., 
\[ 
 \lorb{\sigma} = 
  \Om^{C_{r}}-\Om^{\mathcal{O}}_{(L_0)_M^G}
\]
where $\Om^{\mathcal{O}}_{(L_0)_M^G}$ is the canonically reduced
symplectic form on $(\orb\spr{0}M)_{(L_0)_M^G}$.  

From the general theory we know that the Hamiltonian $\hfree$ reduces
to a Hamiltonian $\hcml$ on the stratum $(T^ {*}V\sporb G)_{(L)}$, and
that integral curves of $\hfree$ project to integral curves of
$\hcml$. 
In particular the dynamics remain confined to the symplectic stratum. 
The reduced Hamiltonian is thus given by 
\[
 \hcml(q,p,[Z]) = \hfree(q,p+A_{q}^{*}(\lam))
\]
where $[Z]$ is the class of $Z$ in $(\orb\spr{0}M)_{(L_0)_M^G}$ and
$A_{q}^{*}: \gu_{q}^{\bot}=\mo^{\bot}\to T_{q}(G.q)=\Sigma^{\bot}$ 
is the point wise dual to the mechanical connection as
introduced in Section \ref{sec:gauged_red}.
Let $q\in C_{r}$, $p=\sum_{i=1}^{l}p_{i}B_{0}^{i}$, and 
$Z=\sum_{\lam\in
  R}\sum_{i=1}^{k_{\lam}}z_{\lam}^{i}E_{\lam}^{i} \in
(\orb\cap\mo^{\bot})_{(L_0)_M^G}$ where $l=\dim\Sigma$ and
$k_{\lam}=\by{1}{2}\dim\lo_{\lam}$. The notation here is as in 
Section~\ref{sec:pol_rep}, 
and $R=R(\lo,\Sigma)\subeq\Sigma^{*}$ denotes the
set of restricted roots, in particular. With these definitions the
dual mapping to the mechanical connection is given by
\[
 A_{q}^{*}(Z)
 = 
 \sum_{\lam\in R}\sum_{i=1}^{k_{\lam}}
  \by{z_{\lam}^{i}}{\lam(q)}B_{\lam}^{i}.
\]
Note that $\lam(q)\neq0$ for all $\lam\in R$ since $q\in C_{r}$ is
regular.
The reduced Hamiltonian thus computes to 
\[
 \hcml(q,p,[Z])
 =
 \by{1}{2}\vv<p+A_{q}^{*}(Z),p+A_{q}^{*}(Z)>
 =
 \by{1}{2}\sum_{i=1}^{l}p_{i}^{2}
  + 
 \by{1}{2}\sum_{\lam\in R}
  \by{\sum_{i=1}^{k_{\lam}}z_{\lam}^{i}z_{\lam}^{i}}{\lam(q)^{2}}.
\]
The reduced Hamiltonian system $(T^{*}V\sporb
G,\sigma^{\mathcal{O}},\hcm)$ is thus a new version of a 
Calogero-Moser system with
spin. It is, in fact, a stratified Hamiltonian system in the sense that
it is a Hamiltonian system on each symplectic stratum $(T^{*}V\sporb
G)_{(L)}$, and the dynamics remain confined to these strata. 

We now show that the thus obtained Calogero-Moser system is
integrable in the non-commutative sense. 
To do so we will use Theorem \ref{thm:zung}. We start by choosing
coordinates $q_{1},\dots,q_{n},p_{1},\dots,p_{n}$ on $T^{*}V=V\times
V$ such that the Poisson bracket of functions $f,g\in\cinf(V\times V)$
is given by the usual equation 
$\bra{f,g} 
 = 
 \sum_{i=1}^{n}
  (\by{\del f}{\del p_{i}}\by{\del g}{\del q_{i}} 
   - \by{\del f}{\del q_{i}}\by{\del g}{\del q_{i}})$.
Moreover we assume that $q_{1}\dots,q_{l},p_{1},\dots,p_{l}$ are
coordinates on $\Sigma\times\Sigma\hookto V\times V$. Let us now
consider the map 
\[
  \Phi: V\times V\longto \Sigma^{\bot}\times V
\]
given by projection, and endow $\Sigma^{\bot}\times V$ with the
inherited Poisson structure. Clearly, $\cinf(\Sigma^{\bot}\times V)$
has a center and this is just generated by $p_{1},\dots,p_{l}$. Thus
we may identify $Z(\cinf(\Sigma^{\bot}\times V)) = \cinf(\Sigma)$. 
Now the set of all first integrals of $\hfree$,
i.e.,
\[
 \FO_{H_{\textup{free}}} 
 =
 \set{F\in\cinf(V\times V): \bra{F,H}=0}
\]
can be identified with $\cinf(\Sigma^{\bot}\times V)$ via
$\Phi$ since $\hfree$ factors over the projection onto the second
factor and is $G$-invariant, and thus can be considered as a function
on $\Sigma$.  
Therefore,
\begin{align*}
 \dim V\times V
 =
 \dim\Sigma^{\bot}\times V + \dim\Sigma
 =
 \ddim\FO_{H_{\textup{free}}}+\ddim Z(\cinf(\Sigma^{\bot}\times V)),
\end{align*}
and we are exactly in the situation of the following theorem to
conclude non-commutative integrability of the reduced system.

\begin{theorem}\label{thm:zung}
Assume the Hamiltonian system $(M,\om,H)$ is invariant under a
Hamiltonian action of a compact Lie group $G$. 
If $(M,\om,H)$ is non-commutatively integrable 
(Definition~\ref{def:int}) then the
reduced system is integrable as well:
\begin{itemize}
\item The singular Poisson reduced system 
  is non-commutatively integrable.
\item The singular symplectic reduced system is non-commutatively
  integrable.
\end{itemize}
\end{theorem}

\begin{proof}
This theorem is proved by Zung \cite[Theorem 2.3]{Zung02}. For
material on singular reduction we refer to Ortega and Ratiu
\cite{OR04} and Section~\ref{s:sg_com_red}.
\end{proof}

The idea of non-commutative integrability under the name of degenerate
integrability is due to Nehoro\v{s}ev \cite{Neh72} 
who also introduced the appropriate concept of action-angle variables. 
This section follows mainly the approach of Zung
\cite{Zung02,Zung03}. See also Mishchenko and Fomenko \cite{MF78}. The
following definition is less general than that given in the above
cited references but better suited for the applications in this paper.

\begin{definition}\label{def:int}
Let $(M,\bra{\phantom{f},\phantom{f}})$ 
be a Poisson manifold, and consider a Hamiltonian
function $H: M\to\R$. 
We denote the Poisson sub-algebra of all first integrals of $H$ by
$\FO_{H}$, that is
\[
 \FO_{H}
 :=
 \set{F\in\cinf(M): \bra{F,H}=0}.
\]
The Hamiltonian system is called
\caps{non-commutatively integrable} if there is a finite dimensional
Poisson vector space $W$ and a generalized \momap $\Phi: M\to W$ which is a Poisson
morphism with respect to the Poisson structure on $W$
such that the following are satisfied.
\begin{itemize}
\item $\Phi^{*}: \cinf(W)\to\FO_{H}$ is an \iso of Lie-Poisson algebras. 
\item $\dim M = \ddim\cinf(W)+\ddim Z(\cinf(W))$ where
  $Z(\cinf(W))$ denotes the commutative sub-algebra of Casimir
  functions on $W$, and $\ddim\cinf(W)=\dim W$ is the functional
  dimension of $\cinf(W)$.  
\end{itemize}
\end{definition}

It is crucial in the formulation of the above theorem that $\dim M =
\ddim\FO_{H} + \ddim Z(\FO_{H})$, and $\FO_{H}$ is the set of
\emph{all} first integrals of $H$.

\section{Appendix: Polar representations}\label{sec:pol_rep}

Let $V$ be a real Euclidean vector space, and $G$ be a connected compact Lie
group. Further, let $\rho: G\to\SO(V,\vv<\_\,,\_>)$ be a \caps{polar
representation} of $G$ on $V$. That is, there is subspace
$\Sigma\subeq V$ (a \caps{section}) such that $\Sigma$ meets all
$G$-orbits, and does so orthogonally. 
The following is due to Dadok \cite{Dad85} and is a consequence of his
classification of polar actions.

\begin{proposition}
There exists a connected Lie group $\tilde{G}$ together with a
representation $\tilde{\rho}:\tilde{G}\to\SO(V)$ such that the
following hold.
There is a real semi-simple Lie algebra $\lo$ with a Cartan
decomposition $\lo = \gu\oplus\po$. Moreover, there is a Lie algebra
isomorphism $A: \lie{\tilde{G}}=\tilde{\gu}\to\gu$ and a linear
isomorphism 
$B: V\to\po$ such that
$B(\tilde{\rho}^{\prime}(X).v) = [A(X),B(v)]$ for all $X\in\tilde{\gu}$ and $v\in
V$. Finally, the $G$-orbits coincide with the $\tilde{G}$-orbits, that
is $V/G=V/\tilde{G}$. 
\end{proposition}

\begin{proof}
See Dadok \cite[Proposition 6]{Dad85}.
\end{proof}

Thus, for the purpose of this paper, it suffices to assume that
the representation of $G$ on $V$ is a symmetric space representation
whence 
$\lo=\gu\oplus V$ is a Cartan decomposition, and hence
$[\gu,\gu]\subeq\gu$, $[\gu,V]\subeq V$, and $[V,V]\subeq\gu$. 
Therefore, $G\times V\cong L$, $(g,v)\mapsto
g\exp(v)$ is a global Cartan decomposition with compact $G$ where $\lie{L}=\lo$. 

An element $v\in V$ is said to be \caps{regular} (with respect to the
$G$-action) if the orbit $\orb(v)=\rho(G).v=G.v$ is of maximal
possible dimension. The set of regular elements will be denoted by
$V_{r}$. The following assertions which are easy to verify are used in 
Subsection~\ref{sub:pol-cm}. 
(See also Knapp \cite[Chapter VI]{Knapp96}.)

Let $v\in V$. Then, by reason of dimension,
$\ad(v)|Z_{\gu}(v)^{\bot}: Z_{\gu}(v)^{\bot}\to Z_{V}(v)^{\bot}$ 
and $\ad(v)|Z_{V}(v)^{\bot}: Z_{V}(v)^{\bot}\to Z_{\gu}(v)^{\bot}$  
both are linear isomorphisms. 

The set $V_{r}$ of regular elements is open dense in $V$. Moreover,
$v\in V_{r}$ if and only if $Z_{V}(v)=:\Sigma$ is a section in $V$. This is the
case if and only if $\Sigma$ is maximally Abelian.

Let $\Sigma\subset V$ be a section, and put $\mo:=Z_{\gu}(\Sigma)$. The
set $R=R(\lo,\Sigma)\subeq\Sigma^{*}$ shall denote the set of restricted
roots. This gives rise to the restricted root space decomposition 
\[
 \lo
= \mo\oplus\Sigma\oplus\oplus_{\lambda\in R}\lo_{\lambda}. 
\]
Any
\csa $\ho\subeq\lo$ of $\lo$ is of the form 
$
 \ho=\too\oplus\Sigma
$
where $\too\subeq\mo$ is a \csa (Lie algebra to a maximal torus) of $\gu$. 

Each restricted root space
$\lo_{\lambda}$ has an orthonormal basis
$
 E_{\lambda}^{i}\in\gu, 
 B_{\lambda}^{i}\in V
$
where  $i=1,\ldots,k_{\lambda}=\by{1}{2}\dim\lo_{\lambda}$, and which is such
that 
$\ad(v)E_{\lambda}^{i} = \lambda(v)B_{\lambda}^{i}$ and 
$\ad(v)B_{\lambda}^{i} = \lambda(v)E_{\lambda}^{i}$ for all 
$v\in\Sigma$. 
The vectors  
$
 E_{0}^{i}, 
$
where $i=1,\dots,\dim\mo$,
and
$ 
 B_{0}^{j}, 
$
where
$
j=1,\dots,\dim\Sigma
$
will denote an orthonormal basis of $\mo$ and $\Sigma$ respectively. 

\end{document}